# Machine Learning Powered Feasible Path Framework with Adaptive Sampling for Black-box Optimization

Zixuan Zhang[1,2], Xiaowei Song[1,2], Jiaming Li[1,2], Yujiao Zeng[1], Yaling Nie[1], Min Zhu[1], Dongyun Lu[1], Yibo Zhang[1], Xin Xiao[1,4,1], Jie Li[3,4,§]

[1] Institute of Process Engineering, Chinese Academy of Science, Beijing, 100190, China
[2] School of Chemical Engineering, University of Chinese Academy of Science, Beijing 100049, China
[3] Centre for Process Integration, Department of Chemical Engineering, School of Engineering, The University of Manchester, Manchester M13 9PL, United Kingdom
[4] Carbon Neutral Intelligent Industry Innovation Working Committee, China Institute for Innovation and Development Strategy, Beijing, 100044, China

## Abstract

Black-box optimization (BBO) involves functions that are unknown, inexact and/or expensive-to-evaluate. Existing BBO algorithms face several challenges, including high computational cost from extensive evaluations, difficulty in handling complex constraints, lacking theoretical convergence guarantees and/or instability due to large solution quality variation. In this work, a machine learning-powered feasible path optimization framework (MLFP) is proposed for general BBO problems including complex constraints. An adaptive sampling strategy is first proposed to explore optimal regions and pre-filter potentially infeasible points to reduce evaluations. Machine learning algorithms are leveraged to develop surrogates of black-boxes. The feasible path algorithm is employed to accelerate theoretical convergence by updating independent variables rather than all. Computational studies demonstrate MLFP can rapidly and robustly converge around the KKT point, even training surrogates with small datasets. MLFP is superior to the state-of-the-art BBO algorithms, as it stably obtains the same or better solutions with fewer evaluations for benchmark examples.

---

[1] Corresponding author: Prof. Xin Xiao, E-mail: xxiao@ipe.ac.cn. Tel: +86 139 1119 6814.
[§] Corresponding author: Dr. Jie Li, E-mail: jie.li-2@manchester.ac.uk. Tel: +44 161 529 3084.

## 1. Introduction

Black-box optimization (BBO) problems are a category of optimization problems where the algebraic form of the objective function and/or constraints is either unknown, inexact, or expensive to evaluate[1–3]. In chemical engineering, many systems cannot be fully captured by mechanistic models such as complex chemical[4–8] or biological[9,10] reaction networks (e.g., multi-step catalytic reaction [8,10]) and molecular structure–activity relationships[11–13] (SARs, e.g., linking polymer chain architecture to material performance[11]). Furthermore, even when mechanistic models is available to describe specific processes, solving them is frequently prohibitively expensive, which is true of molecular simulations (e.g., density functional theory calculations)[14–16], flow field simulations (e.g., computational fluid dynamics for reactor design) [17–19], and process simulations (e.g., rigorous distillation column modeling)[20–22]. Such systems are thus treated as black-boxes, yet they serve as critical foundational components in chemical design and manufacturing. Optimizing their design and operation typically depends on optimization approaches[4–6,17,23–25]. Therefore, advancing BBO methods can substantially accelerate research and development of chemical engineering.

However, BBO faces a series of challenges. Owing to the absence of closed-form expressions for black-boxes, it is quite difficult to evaluate gradient information of the black-box functions for generating an effective search direction. As a result, researchers usually evaluate the objective function value directly and determine the next iteration point based on the current function value rather than its derivatives[26–28]. A well-known category of BBO algorithms is stochastic or heuristic algorithms, such as hit-and-run algorithm[29], particle swarm optimization (PSO)[30,31], simulated annealing (SA)[28] and genetic algorithm (GA)[32,33]. These algorithms search and break away from a local optimum by designing random or heuristic rules. Although they provide global search capability of finding good feasible solutions, convergence cannot be theoretically guaranteed. They also need extensive evaluation of the black-box functions, substantially increasing the computational effort.

Correspondingly, some BBO algorithms employ deterministic partition search rules, such as pattern search[34], multilevel coordinate search[35], DIRECT (DIviding RECTangles) algorithm[26,36]. However, these algorithms typically require extensive evaluations of black-box functions, resulting in the optimization extremely time-consuming or expensive. Therefore, researchers have adopted surrogate model-based methods to address BBO problems, reducing the number of evaluations of black-box functions.

Surrogate model-based BBO methods can be broadly categorized into two classes: one involves sampling and updating the surrogate models dynamically during the optimization process; the other relies on constructing a high-fidelity surrogate model in advance and then performing optimization directly based on it. The former mainly includes implicit filtering[37], rectified linear unit neuron network (ReLU NN)-based optimization approach with adaptive uncertainty-aware sampling[38], the automated learning of algebraic models for optimization (ALAMO)[39,40], trust-region methods[41–44], and Bayesian optimization (BO)[45–49]. The implicit filtering algorithm is designed to solve bound-constrained BBO problems by constructing the approximated derivatives through numerical differentiation[50]. However, this method is not well-suited for general BBO problems. The ReLU NN-based optimization framework[38] initially employs low-discrepancy sequence sampling and trains a ReLU NN as surrogate. At each iteration, adaptive resampling with uncertainty assessment is conducted to update the ReLU NN, and the problem is reformulated as MILP for solution. Complex neural network architectures can lead to an excessive number of binary variables, making the problem computationally intractable or even infeasible[51]. This implies that the structure of the neural networks should remain relatively simple, making it suitable only for low-dimensional or weakly nonlinear problems. ALAMO employs best subset selection techniques to construct low-complexity algebraic models, which are iteratively updated through adaptive sampling during the optimization process. ALAMO needs to solve a mixed-integer linear programming (MILP) problem when constructing the surrogate model at each iteration. The accuracy of the surrogate model is critically dependent on the pre-specified basis functions. More importantly, the subset selection methods may not perform well under high-dimensional inputs.[39]

The trust-region methods offer a theoretically grounded approach for solving general BBO problems. When employing κ-fully linear models, trust-region filtering (TRF) has been proven to converge to first-order critical points[42,52]. Nevertheless, the method could require a substantial number of function evaluations in practice[43]. To address this issue, Liang et al. proposed a TRF variant based on Gaussian Process (GP) surrogate models[43], which was able to converge to regions close to the KKT points using significantly fewer evaluations. One limitation of the TRF methods is that, when analytical gradients are unavailable, it resorts to finite differences for gradient estimation[53], this can introduce inaccuracies in the Jacobian matrix. In turn, the approximate Hessian matrix used in TRF may become ill-conditioned, undermining stability and convergence precision. Besides, TRF may exhibit limited applicability in scenarios involving multiple black-box functions, as it would require establishing separate trust regions and sampling subspace for each black-box function.

BO is another widely used[45-49] BBO algorithm due to its capability to quantify the uncertainty of black-box functions as well as to establish surrogate models by using only a small amount of samples. The most commonly used surrogate model is Gaussian processes (GP)[54] in BO where the noise-free GP model is also known as Kriging model[55,56]. As a non-parametric model, GP is computationally efficient under a small dataset, making it well-suited for situations where frequent updates of surrogate models are required during the optimization process. Most BO-based algorithms assume only inequality constraints are involved in the generally constrained BBO problems and have little emphasis on BBO problems with complex equality constraints. Recently, Tian and Ierapetritou[55] proposed a framework for BBO problems with inequality constraints, which are formulated as a classification task. Based on this framework, Fan et al.[53] further developed an adaptive

**4**

sampling Bayesian optimization algorithm (ASBO) to handle both inequality and equality constraints. Additional strategies were introduced to achieve a balanced trade-off between solution feasibility and optimality. However, a notable limitation of BO lies in the stochasticity of its outcomes. Computational studies have shown that, with different initial data, repeated runs of BO can yield results that vary significantly — in some cases, by up to 300%.[53]

To mitigate the sensitivity of GP-based method (including BO and GP-based trust-region method) to initial sampling, the alternative method that first construct a relatively accurate surrogate model and then execute the optimization would be preferred. This alternative method is especially effective in certain situations where the available dataset is quite limited, and active sampling may not be viable in practice. Under such a framework, researchers focus on constructing more accurate surrogate models and attempt to find the theoretically optimal solutions of surrogates. With the advancement of machine learning (ML) techniques, taking ML models as surrogates and then optimizing has been widely adopted by the PSE community due to their fast output prediction and strong regression capabilities. Therefore, algorithms that rely on extensive function evaluations to find an optimum, such as GA, and PSO, are commonly used to optimize the surrogates[57–59]. More recently, researchers have made significant effort in combining ML surrogates with deterministic derivative-based algorithms[23,60–62], offering a theoretical convergence. Among the most advanced techniques, researchers reformulate the explicit algebraic expressions of ML models in conjunction with deterministic solvers[62–68]. Some researchers directly employ the explicit formulations of neural networks as nonlinear constraints[65–67]. Concurrently, another cohort of researchers has explored the utilization of decision trees [61–63] or neural networks with ReLU activation functions [51,64,69–71], which can be reformulated as mixed integer constraints. This strategy aims at tackling the ML surrogates assisted BBO with MILP[63,70] or transforming nonconvex NLP into MILP which can be solved to global optimality[51,62]. However, these methods usually neglect the computational burden for sampling and training the ML models, which is usually time-consuming. A trade-off among the number of samples, prediction accuracy, and the complexity of algebraic expressions should be balanced well. This is because higher model accuracy generally requires more samples and introduces more nonlinear terms or binary variables, thereby increasing the challenges for optimization.

To overcome disadvantages of the aforementioned BBO methods, in this study we propose a novel machine learning-based feasible path optimization (MLFP) framework for solving general BBO problems where complex equality and inequality constraints are involved. In this framework, an efficient adaptive sampling strategy is first proposed, which filters out potentially infeasible samples with support vector machine (SVM) and focuses on regions likely to contain optimal solutions. This ensures that an informative set of samples can still be obtained even under limited data conditions. As a result, a small number of samples (e.g., 1000) are required to develop high prediction accuracy of the surrogate models by machine learning algorithms. Then, the feasible path optimization methods[72–75] are integrated, which enables the optimization process to circumvent large number of intermediate variables within the surrogate models, thereby obtaining derivatives of outputs respect to inputs. As a result, even sophisticated ML models can be effectively employed as surrogates without compromising computational costs and predictive accuracy. The computational results demonstrate that the proposed adaptive sampling approach achieves a substantial reduction in the required number of samples compared to widely used Latin hypercube sampling (LHS) strategy. The implementation of the feasible path method significantly improves computational speed, allowing the optimization to be completed within just a few seconds. For all benchmark problems with known global optima, the proposed MLFP reliably converges to the true optimal solution, highlighting its effectiveness and numerical robustness. It is also superior to the state-of-the-art BBO algorithms due to the same or better solutions being stably obtained with fewer evaluations for the benchmark examples.

The remaining structure of the article is as follows. The overall NLP BBO problem studied is defined in Section 2. The detailed information of surrogate-assisted hybrid feasible path framework is presented in Section 3. Section 4 offers thirteen Examples and their results as well as comparison with state-of-the-art BBO algorithms. Finally, Section 5 provides conclusion and prospect.

## 2. Problem description

The general form of a BBO problem to be solved in this work is presented below, which is a nonlinear programming (NLP) problem.

$$
\begin{aligned}
&\min_{\mathbf{x}\in\mathbb{R}^n} \quad f(\mathbf{x}) \\
&\text{s.t.} \quad \mathbf{h}(\mathbf{x}) = 0 \\
&\qquad\;\; \mathbf{g}(\mathbf{x}) \leq 0 \\
&\qquad\;\; \mathbf{x}^{lb} \leq \mathbf{x} \leq \mathbf{x}^{ub}
\end{aligned}
\tag{P0}
$$

where the function $f : \mathbb{R}^n \to \mathbb{R}$ is the objective function, $\mathbf{h} : \mathbb{R}^n \to \mathbb{R}^{p_E}$ is the equality constraints and $\mathbf{g} : \mathbb{R}^n \to \mathbb{R}^{p_I}$ is the inequality constraints excluding the bounding constraints. $p_E$, and $p_I$ represent the dimensions of equations and inequalities respectively; $\mathbf{x}$ represent a vector of variables; the superscript $ub$ represents the upper bound, and $lb$ represents the lower bound. In this setting, the black-box function may pertain to the objective function $f$ as well as to the constraint functions $\mathbf{h}$ (equalities) and $\mathbf{g}$ (inequalities), depending on the problem structure.

Without loss of generality, we divide the equality constraints $\mathbf{h}(\mathbf{x}) = 0$ into $\mathbf{h_1}(\mathbf{x}) = 0$ and $\mathbf{h_2}(\mathbf{x}) = 0$ where $\mathbf{h_1}(\mathbf{x}) = 0$ containing all black-box constraints and $\mathbf{h_2}(\mathbf{x}) = 0$ without containing any black-box functions. If any inequality constraints $\mathbf{g_j}(\mathbf{x}) \leq 0$ containing black-box functions, we can introduce auxiliary variables $\mathbf{v}$ and let $\mathbf{v} = \mathbf{g_j}(\mathbf{x})$. We then add $\mathbf{v} = \mathbf{g_j}(\mathbf{x})$ into $\mathbf{h_1}(\mathbf{x}) = 0$ and remove $\mathbf{g_j}(\mathbf{x}) \leq 0$ from the inequality constraint sets $\mathbf{g}(\mathbf{x}) \leq 0$. The auxiliary variable $\mathbf{v}$ is then added into variables $\mathbf{x}$. The original BBO problem (**P0**) can be reformatted into (**P1**) below,

$$\min_{\mathbf{x}\in\mathbb{R}^n} \quad f(\mathbf{x})$$
$$\text{s.t.} \quad \mathbf{h_1}(\mathbf{x}) = 0$$
$$\mathbf{h_2}(\mathbf{x}) = 0 \qquad\qquad \textbf{(P1)}$$
$$\mathbf{g}(\mathbf{x}) \leq 0$$
$$\mathbf{x}^{lb} \leq \mathbf{x} \leq \mathbf{x}^{ub}$$

where $\mathbf{h_1}\colon \mathbb{R}^n \to \mathbb{R}^{p_{E1}}$ and $\mathbf{h_2}\colon \mathbb{R}^n \to \mathbb{R}^{p_{E2}}$ are the equality constraints and $\mathbf{g}\colon \mathbb{R}^n \to \mathbb{R}^{p_I}$ is the inequality constraints excluding the bounding constraints. $p_{E1}$, $p_{E2}$ and $p_I$ represent the dimensions of equations and inequalities respectively. In this setting, the black-box function may pertain to the objective function $f$ as well as to the constraint functions $\mathbf{h_1}$ (equalities).

## 3. Surrogate-assisted feasible path framework

As $f(\mathbf{x})$ and/or $\mathbf{h_1}(\mathbf{x}) = 0$ in problem (**P1**) contain black-box functions that are unknown or expensive for evaluation, we need to develop surrogate models for these black-box functions. To do this, we can divide all variables in (**P1**) into independent variables (denoted as $\mathbf{x}_I$) and dependent variables (denoted as $\mathbf{x}_D$) where $\mathbf{x}_D$ has the same dimension as $\mathbf{h_1}$. Then, problem (**P1**) can be transformed into problem (**P2**).

$$\min_{\mathbf{x}\in\mathbb{R}^n} \quad f(\mathbf{x}_I, \mathbf{x}_D)$$
$$\text{s.t.} \quad \mathbf{h_1}(\mathbf{x}_I, \mathbf{x}_D) = 0$$
$$\mathbf{h_2}(\mathbf{x}_I, \mathbf{x}_D) = 0 \qquad\qquad \textbf{(P2)}$$
$$\mathbf{g}(\mathbf{x}_I, \mathbf{x}_D) \leq 0$$
$$\mathbf{x}_I^{lb} \leq \mathbf{x}_I \leq \mathbf{x}_I^{ub}$$
$$\mathbf{x}_D^{lb} \leq \mathbf{x}_D \leq \mathbf{x}_D^{ub}$$

The proposed optimization framework is illustrated in FIGURE 1. The proposed framework mainly consists of three components, including the adaptive sampling strategy, construction of surrogate models, and the feasible path-based optimization approach. The details about each component will be explained in following subsections.

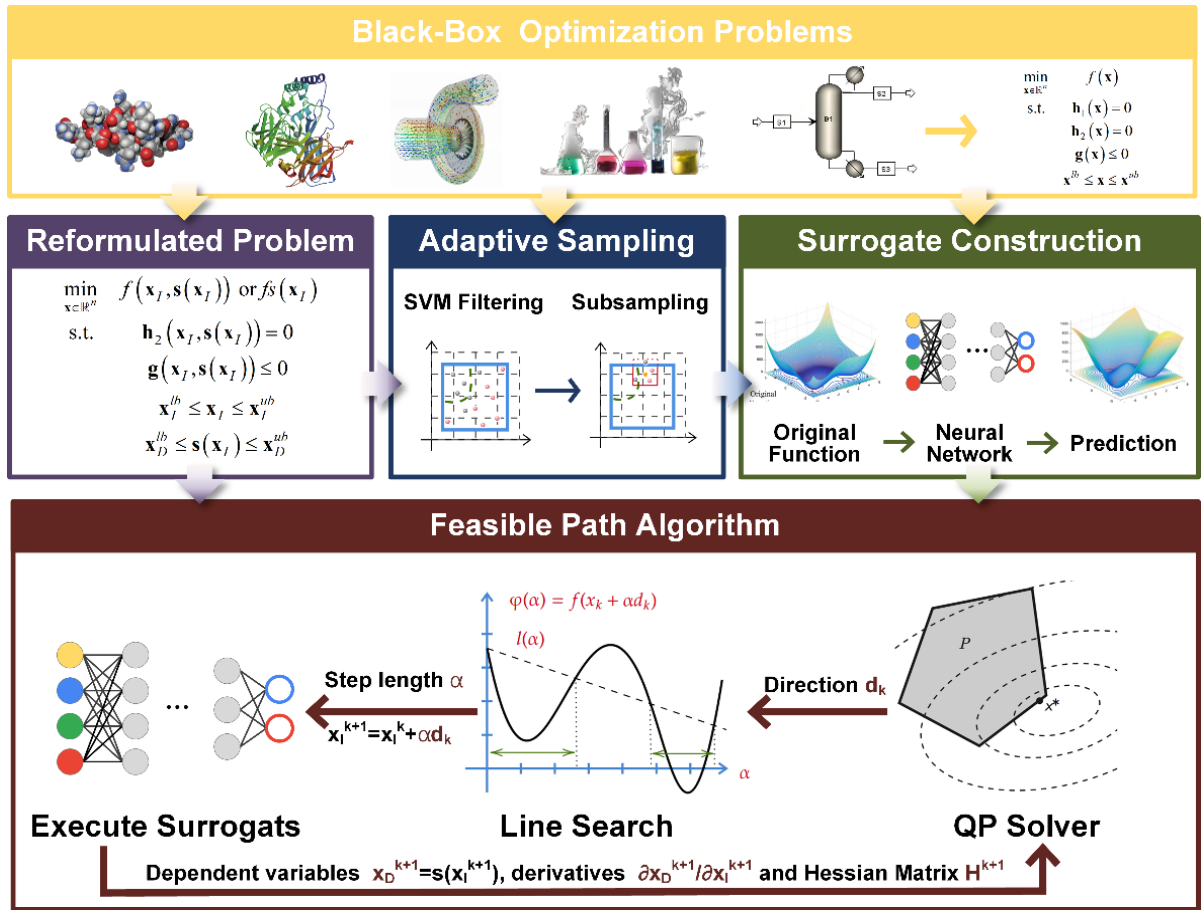

FIGURE 1 Proposed MLFP framework for general BBO problems.

## 3.1 Adaptive sampling strategy

Adaptive sampling plays an increasingly important role in surrogate assisted BBO[9,23,53], due to its ability to generate higher quality data, establish more accurate models, and reduce the number of function evaluations. Most adaptive sampling strategies are typically integrated with the optimization process iteratively[23,39,53,76], whereby the next sampling point can be selected based on the current optimization results. While effective for algorithms coupling optimization with sampling, such approaches are not well-suited for applications in the decoupled algorithms. To address this limitation, we propose a new adaptive sampling strategy to construct a global dataset that can be used to establish surrogate models before optimization, with particular emphasis on capturing fine details near the optimal region and bounds.

In the proposed strategy, we consider two scenarios. The first scenario is the case where the outputs of the surrogate model constructed are not constrained in the optimisation model. In other words, all dependent variables in the surrogate models are free variables in the optimisation problem. The other scenario is the case where the outputs of the surrogate model established are constrained in the optimisation problem. While we focus only on the potentially optimal region and the details near the bounds for the first scenario, we must further consider the regions close to the constraint boundaries for the latter. The existence of constraints on the outputs of the surrogate models is clearly a more general consideration, whereas the first scenario can be handled with a simplified approach. We provide a detailed description of the adaptive sampling strategy proposed in this work in the following.

Firstly, the original bounds of each independent variable are slightly expanded near its boundary to form new bounds (e.g., $\pm$ 5% of the original bounds). Then sampling on the independent variables is conducted within the new bounds by using LHS and a small-scale sample of independent variables (e.g., 10% of the total samples required) is thus obtained. At each sample, the black-box function is evaluated with the purpose of obtaining an overall understanding of the black-box function values and identifying potentially promising regions.

Secondly, the samples generated are used to evaluate constraints to check their feasibility. If a sample satisfies all constraints, then it is labelled as 1, indicating a feasible point. Otherwise, it is labelled as 0, indicating an infeasible point. All these samples are used as training points to develop a SVM for classification. Due to the limited sample size, the initial SVM classifier tends to be inaccurate. To address this issue, we extract the support vectors from the trained SVM model and retrain the classifier by treating all support vectors as feasible. This approach enables inclusion of a small number of infeasible points near the constraint boundaries within the feasible region, thus ensuring a more accurate representation of the boundary details.

Thirdly, a sub-region is established around the best-performing point currently identified within the feasible region, and LHS is conducted within this localized area. Importantly, we employ the trained SVM model to exclude potentially infeasible points from the sample pool, preventing the waste of computational resources on evaluating function values at non-promising, infeasible locations. The computational procedure of the adaptive sampling strategy is illustrated in FIGURE 2. For a comparative illustration against LHS, refer to FIGURE S1 in the **Supplementary Material.**

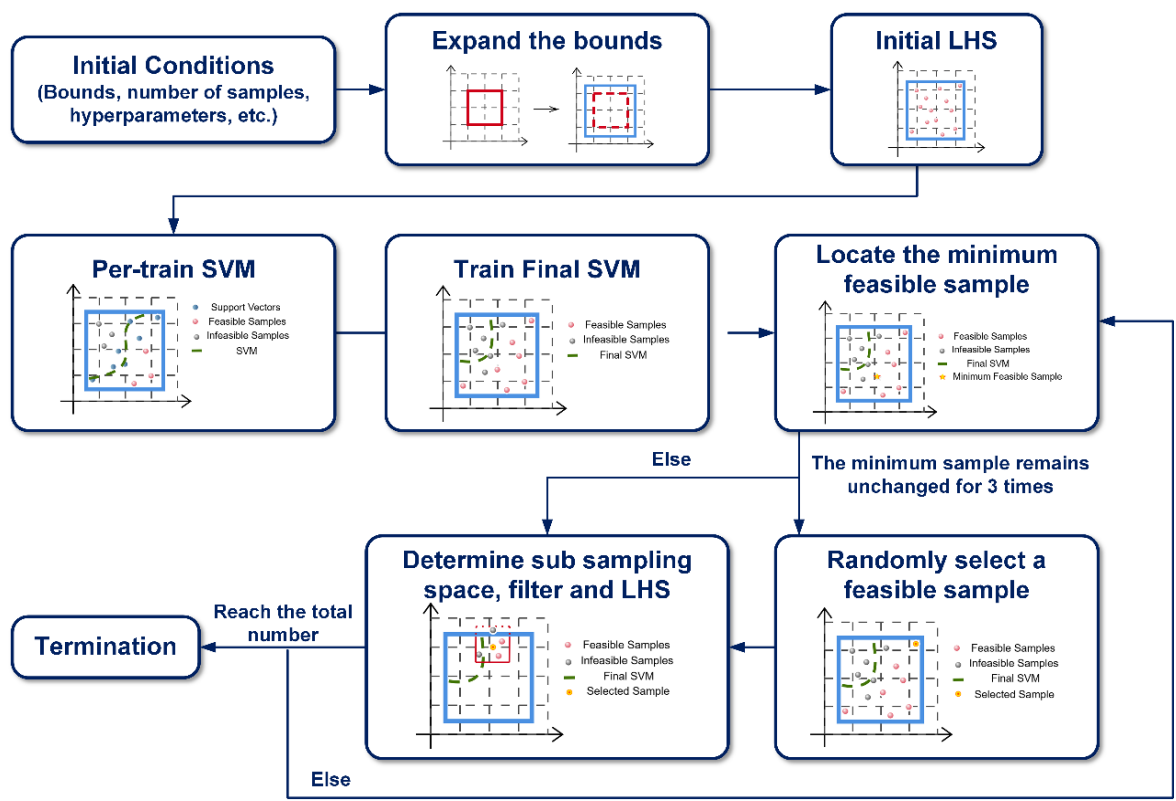

FIGURE 2 Flowchart of the proposed adaptive sampling strategy.

## 3.2 Construction of surrogates with derivatives

Before we build high-quality surrogate models, it is of great importance to determine independent variables and their valid bounds. In this section, we outline the determination of independent variables and their bounds, and explain the criteria for selecting appropriate surrogate models within the proposed framework.

### 3.2.1 Selection of independent variables and bounds

Defining appropriate independent variables and their ranges is crucial for construction of surrogates. Independent variables are typically those that can be manually adjusted, such as reaction temperature, pressure, and residence time. For problems in which the independent variables are not readily apparent, it is essential to prioritize inputs that have a substantial influence on the dependent variables. This can be achieved through degree of freedom analysis, sensitivity analysis and/or feature selection techniques based on a small amount of sampling data. The number of independent variables required for the system can be determined through degrees of freedom analysis. Then sensitivity analysis can be performed to select variables that have significant impacts on the system. Feature selection can also be used to determine the optimal subset of features, including Filter, Wrapper and Embedded methods.

The ranges of independent variables are mainly determined by the following factors. First, there are physical constraints from the real world, including the range of human operability, equipment capacity, design margin, and safe operating conditions. These often define the hard boundaries within which the system can function properly. Second, the reaction conditions themselves impose limitations. For example, a certain chemical reaction may only produce the desired product within a specific temperature range. Third, operational experience or results obtained through trial or experience can also influence variable ranges. While this factor is not mandatory, it is often effective in narrowing down the feasible region and improving optimization efficiency.

### 3.2.2 Selection of appropriate surrogates

Machine learning (ML) models are of particular interest due to their strong regression capabilities. Several commonly used ML models—including decision trees (DT), support vector regression (SVR), neural networks (NN), and Gaussian process regression (GPR)—are introduced. The corresponding algebraic formulations, along with their derivatives, are provided in Section S2 of the **Supplementary Material**.

To ensure mathematical rigor, the surrogate model $s(\mathbf{x})$ should belong to $\mathcal{C}^2$, which means it needs to be both continuous and twice differentiable. Therefore, Decision trees are not recommended for use. This is primarily because their decision boundaries are non-smooth, rendering them non-differentiable at these decision boundaries. Additionally, their predictive models are typically limited to simple forms such as mean-based or linear models, failing to provide second-order derivatives. Besides, given that surrogate models are typically multi-input multi-output systems, NN and GPR are more recommended than SVR. Also, it should be noted that both the kernel function of GPR and the activation function of NN should belong to $\mathcal{C}^2$.

In this work, we choose NN with Swish function as activation function. Compared to GPR, NNs are capable of achieving higher prediction accuracy. This advantage can be attributed to the Universal Approximation Theorem, which theoretically demonstrates that a NN with a single hidden layer can approximate any continuous function. Notably, deep neural networks exhibit greater efficiency in approximating complex, high-dimensional functional relationships, further enhancing prediction performance relative to GPR. In addition, saturated activation functions such as Sigmoid and Tanh tend to suffer from the vanishing gradient problem as network depth increases. To address this issue, a non-saturated activation function belonging to the $\mathcal{C}^2$ is required, which is why we opted for the Swish function.

### 3.3 Feasible path optimization algorithm

After the surrogate models for $\mathbf{h}_1(\mathbf{x}) = 0$ and/or the objective function are developed, then $\mathbf{h}_1(\mathbf{x}) = 0$ and/or the objective function in problem (**P2**) are replaced by the surrogates denoted as $x_D = \mathbf{s}(x_I)$ **and/or** $fs(\mathbf{x}_I)$ and problem (**P3**) below is obtained.

$$
\begin{aligned}
&\min_{\mathbf{x} \in \mathbb{R}^n} && f(\mathbf{x}_I, \mathbf{x}_D) \text{ or } fs(\mathbf{x}_I) \\
&\text{s.t.} && \mathbf{x}_D = \mathbf{s}(x_I) \\
& && \mathbf{h}_2(\mathbf{x}_I, \mathbf{x}_D) = 0 \\
& && \mathbf{g}(\mathbf{x}_I, \mathbf{x}_D) \leq 0 \\
& && \mathbf{x}_I^{lb} \leq \mathbf{x}_I \leq \mathbf{x}_I^{ub} \\
& && \mathbf{x}_D^{lb} \leq \mathbf{x}_D \leq \mathbf{x}_D^{ub}
\end{aligned}
\tag{P3}
$$

To solve problem (**P3**) effectively, we employ the pseudo-transient continuation (PTC) based hybrid feasible path optimization algorithms developed by Ma et al.[73], which divide the entire problem (**P3**) to two subproblems including a small optimization problem in the outer level and a process simulation problem in the inner level. While the small optimization problem in the outer level can be solved using the existing gradient-based optimization algorithms such as improved SQP[72,75], the process simulation problem in the inner level could be solved using solution algorithms for nonlinear equations such as a Newton-based algorithm[77] or PTC approach[78]. More details about the hybrid feasible path optimization algorithms can be referred to Ma et al.[73]. The entire hybrid feasible path optimization algorithms are illustrated in ALGORITHM S2 of the **Supplementary Material**.

## 4. Computational studies

We solve thirteen examples to evaluate the capability of the proposed MLFP framework, including ten benchmark test functions with known globally minimum values[43,79], the classic Williams-Otto process optimization problem[43], extractive distillation for separation of toluene form n-heptane[70] and $CO_2$ capture from biogas from Xu et al.[80]. Throughout these examples, multilayer perceptrons (MLPs) are used to construct surrogate models with the Swish activation function shown in Eq.(1). All examples are solved on a desktop with Intel(R) Core(TM) i9-14900K (3.20 GHz) processor and 128 GB of RAM running Windows 11 64-bit operating system.

$$Swish(x) = \frac{x}{1+e^{-x}} \tag{1}$$

## 4.1  Examples 1-6

Six typical functions from Virtual Library of Simulation Experiments[79] as shown in TABLE S1 of the **Supplementary Material** are selected as test functions. Examples 1-6 are Sphere function, Quadratic function, six hump-camel function, Schaffer No.2 function, Griewank function, and Ackley function, respectively. Except for the Sphere function in Example 1 and a special form of Quadratic function in Example 2, which are high-dimensional convex functions, the other four functions in Examples 3-6 are nonconvex and possess multiple local minimum points. We construct the surrogate models for these 6 functions, which are included in the dataset (the surrogates of the following examples are also included), with the corresponding access link provided in the **Supporting Information**. The performance of the six surrogate models is illustrated in FIGURE S5. From FIGURE S5, it is clearly demonstrated their capability to serve as effective and accurate surrogates for the original test functions (i.e., mean squared error on test set varying from 0.27 to less than 0.01, and $R^2$ varying from 0.87 to 1).

FIGURE 3 shows the convergence curves for solving the six surrogate models by using MLFP with all initial points at $x_i = x_i^{ub}$ or $x_i^{lb}$. $\|\hat{x}^* - x^*\|_2$ represent the distance between the obtained optimal solution $\hat{x}^*$ and true globally optimal solution $x^*$ (marked with black dashed line in FIGURE 3), which does not exceed 0.1 for the six examples. The colored dashed lines represent the predicted values of the surrogate models at each iteration, while the colored solid lines represent the actual output values of the black-box function corresponding to the inputs. According to the true function value at the obtained solution and the globally minimum value, the results indicate that MLFP successfully converges to the global optimum in an extremely short time, varying from 0.1s to 1s in Example 1-6.

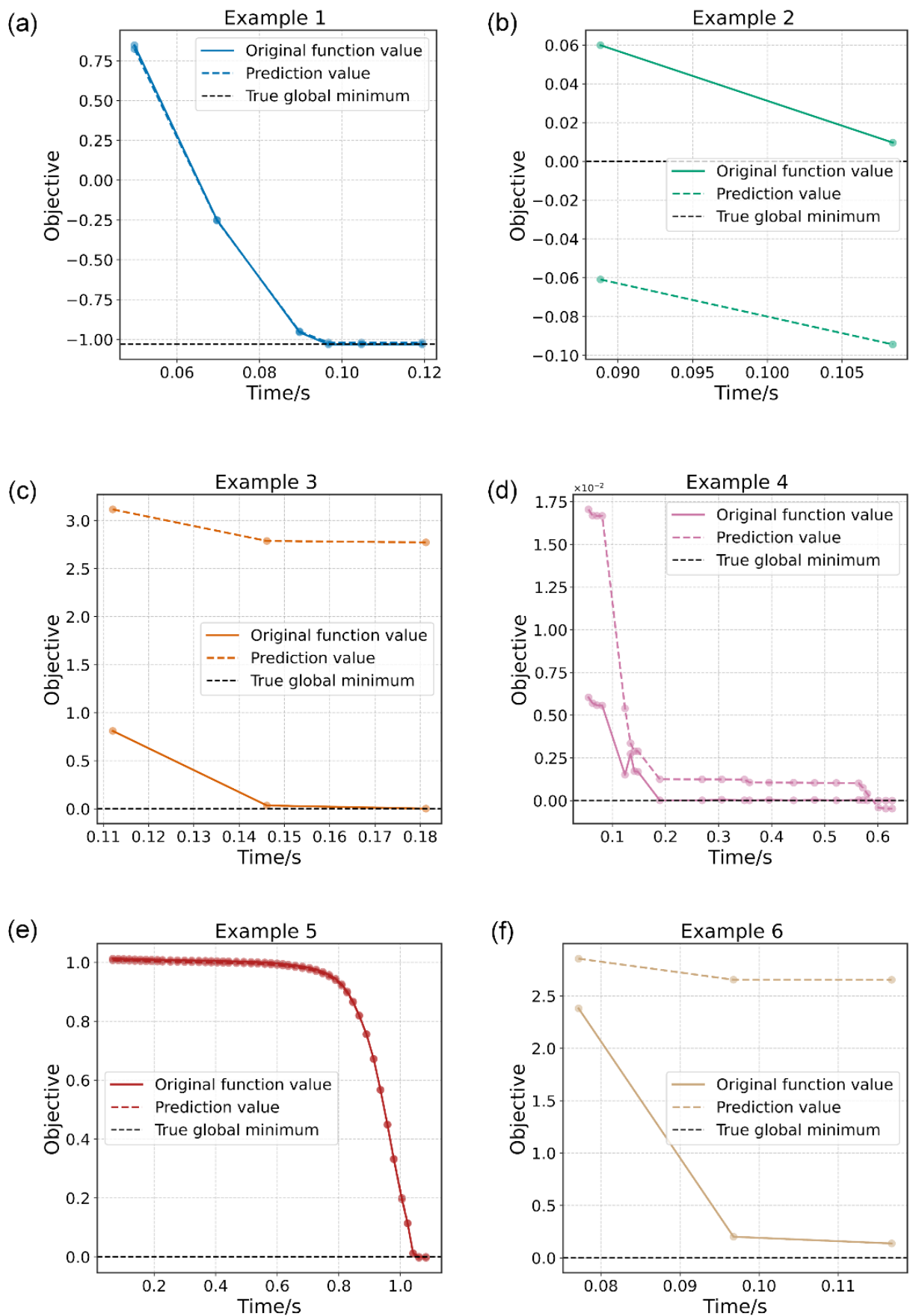

FIGURE 3 Convergence curves of MLFP for Examples 1-6.

## 4.2 Examples 7-10

To further demonstrate the effectiveness of the proposed MLFP framework, we solve 4 additional test functions[43] (i.e., Examples 7-10), including Hartmann 3D function (Example 7), Powel function (Example 8), Rosenbrock function (Example 9) and Trid function (Example 10). We assume that evaluations of these functions were particularly expensive. Therefore, they are treated as black-box functions. We employ the adaptive sampling strategy proposed in the MLFP framework to generate a small yet informative sample set for constructing surrogate models of these black-box functions. The convergence curves are depicted in FIGURE 4 (a-d). As shown in FIGURE 4 (a-d), MLF converges to the optimal solution[79] (see the black dashed line) for all these examples. Specially, MLFP converges to −3.86 for Example 7 within 0.2 s, to 0.08 for Example 8 within 0.35 s, to 0.05 for Example 9 within 0.2 s and to −49.99 for Example 10 within 0.2 s.

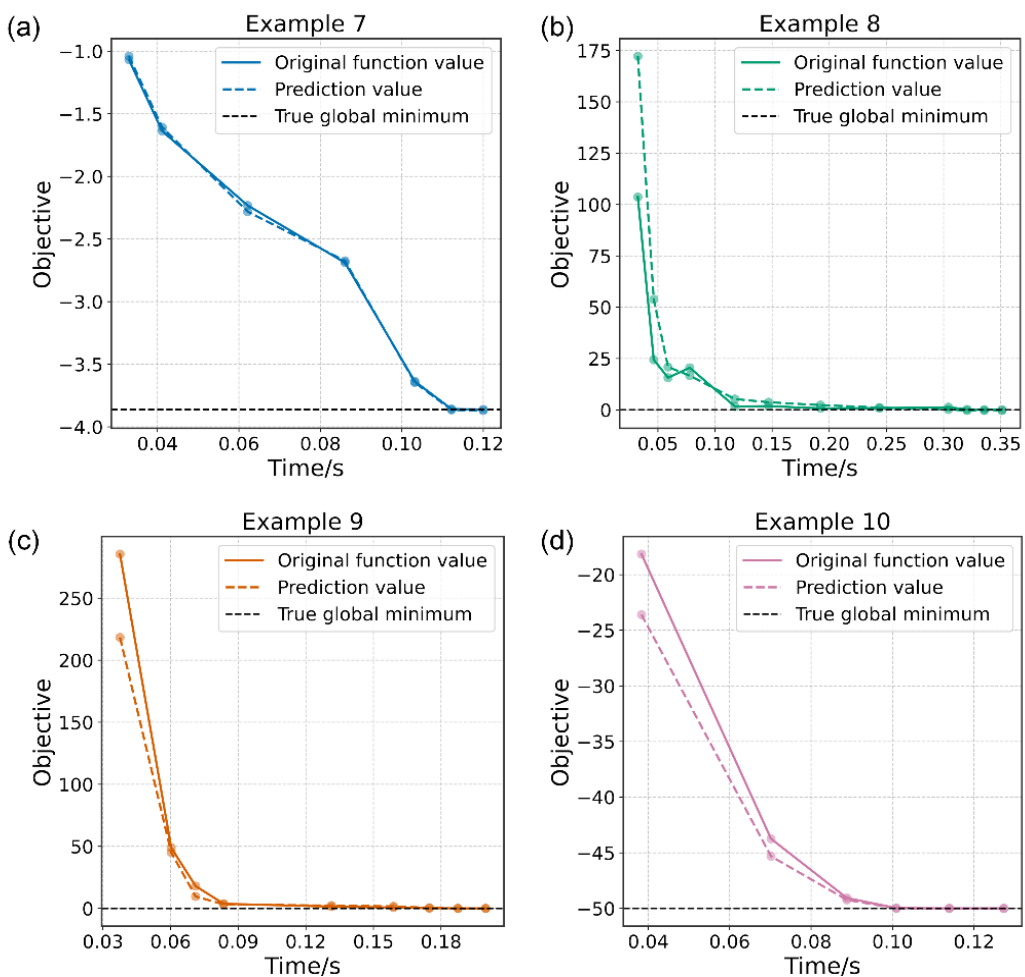

FIGURE 4 Convergence curves of MLFP for Examples 7-10.

Due to the inherent random characteristics of the sample points generated by the proposed adaptive sampling strategy, the surrogate models may be different for each run with varying prediction accuracy, which may affect the optimization performance. To evaluate the stability of the proposed MLFP framework, it is executed 30 times, which are same as that in Liang et al. [43]. The computational results on the number of function evaluations with 30 runs are illustrated in FIGURE 5. Boxplot of the optimization results with 30 runs are illustrated in FIGURE 6. These results are compared with those from TRF-linear (TRF-1), TRF-simplified quadratic (TRF-2), TRF-standard quadratic regression (TRF-3) and TRF-GP (TRF-4) reported in Liang et al. [43].

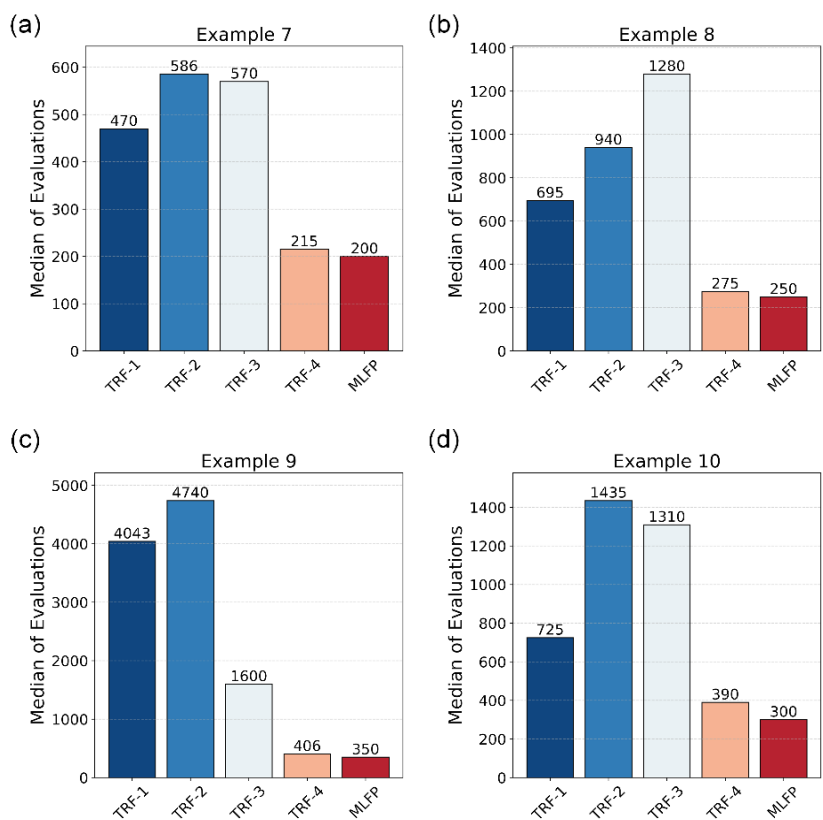

FIGURE 5 Median of function evaluations by TRF 1-4 and MLFP for Examples 7-10.

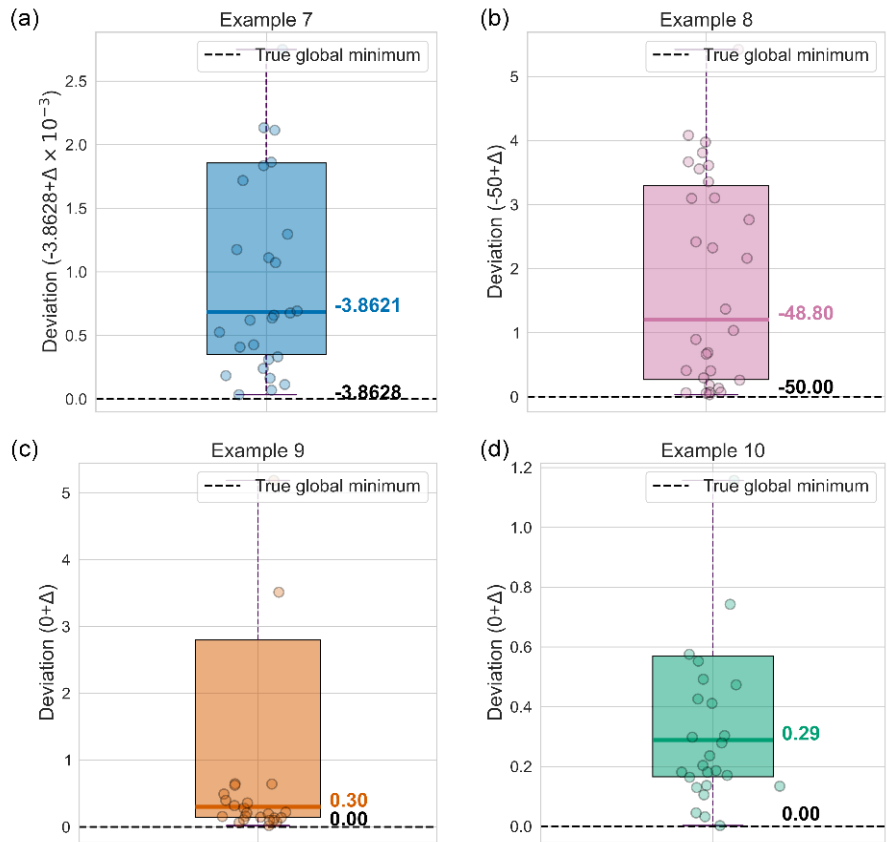

FIGURE 6 Boxplot of optimization results for Examples 7-10 from MLFP with 30 runs.

As illustrated in FIGURE 5(a-d), MLFP requires 200 function evaluations for Example 7, 250 for Example 8, 350 for Example 9 and 300 for Example 10. The number of function evaluations required by MLFP is significantly less than that from TRF-1, TRF-2, and TRF-3. For instance, MLFP requires 350 function evaluations for Example 9, whilst TRF-1, TRF-2 and TRF-3 require 4043, 4740, and 1600 function evaluations, respectively, which are 91.34%, 92.62%, 78.13% more than that from MLFP. The number of function evaluations required by MLFP is still smaller than that required by TRF-4. For instance, MLFP requires 350 function evaluations for Example 9, while TRF-4 requires 406 function evaluations.

From FIGURE 6, it can be observed that MLFP converges to the vicinity of the global optimum for Examples 7 and 10 in most runs, demonstrating its stability. The maximum deviation is about 0.07% for Example 7. However, for Examples 8-9 the median deviation is about 2.4%, $3 \times 10^5$ %, and $2.93 \times 10^5$ %, assuming the minimum value of 0 is $10^{-5}$. The large median deviation for Examples 8-9 is due to the minimum value being 0. Furthermore, as pointed out by Liang et al.[43], for large-scale gray-box problems, the TRF framework requires solving a trust-region subproblem, which is a large-scale NLP, at each iteration. The computational cost of solving such subproblems can even exceed that of training GP. In some cases, the time for a single iteration may even exceed that of training the GP model, potentially taking several minutes per iteration (not including function evaluation time)[43]. In comparison, after the surrogate model has been well-constructed, MLFP only needs to solve the large-scale NLP problem once. The total time of MLFP for Examples 7-10 is within 15 min.

## 4.2 Example 11: Williams-Otto process

The Williams-Otto (WO) process is a classical benchmark that has been widely adopted for evaluating optimization algorithms. A schematic diagram of the WO process and the equation-oriented (EO)-based optimization model can be found in Eq.(S46)-(S52) of the **Supplementary Material**. The objective is to maximize return on investment (ROI) by adjusting the reactor volume $V$, reaction temperature $T$, and the feed flow rates of components A and B ($F_A$ and $F_B$), and stream recycle ratio $\eta$. The inputs to the surrogate models are $V$, $T$, and $F_A$, $F_B$ and $\eta$, which also serve as the independent variables in the optimization problem. The outputs include the objective function value (ROI) and the product stream flow rate $F_P$, which are dependent variables subject to bound constraints.

In this example we compare the performance of TRF 1-4 with MLFP under different samples used. We also compare the performance of the proposed adaptive sampling strategy and the LHS strategy in MLFP. In MLFP with LHS we generate 100, 200, 350, 500, and 1000 valid sampling points (successful simulation convergence), respectively, which correspond to MLFP1-5. In MLFP with the proposed adaptive sampling strategy, we also generate 100, 200, 350, 500 and 1000 valid sampling points, respectively, which are labelled as MLFP6-10.

The convergence curves of MLFPs are shown in FIGURE 7. As shown in FIGURE 7

(a), MLFP7-MLFP10 are capable of converging to the points around the KKT point (see the black dashed line), while MLFP1-5 fails to converge to the KKT point as shown in FIGURE 7 (b). The shaded region around a convergence curve in FIGURE 7 represents the difference between the predicted and actual values. It is clearly shown that when the number of valid sampling points increases, the shaded region becomes narrow, indicating that the difference is reduced. In other words, increasing the number of valid sampling points would help increase model prediction accuracy. It can also be observed that when the iterate approaches the true optimal value, the shaded region diminishes. This may indicate that the samples from the adaptive sampling strategy are mainly clustered around the rigorous-based optimal solution. To verify this hypothesis, we merge the samples obtained from 30 runs, divided the MLFP into 5 groups, and assigned the same number of sampling points to each group. We then conducted sampling for each group using LHS and AS respectively, with the probability density functions shown in FIGURE 8. FIGURE 8 shows that under the same sampling number, the highest peak of AS is closer to 0 than LHS, indicating that the sample distribution is closer to the rigorous model-based optimal solution. And as the number of samples increases, the peak of AS becomes closer to 0.

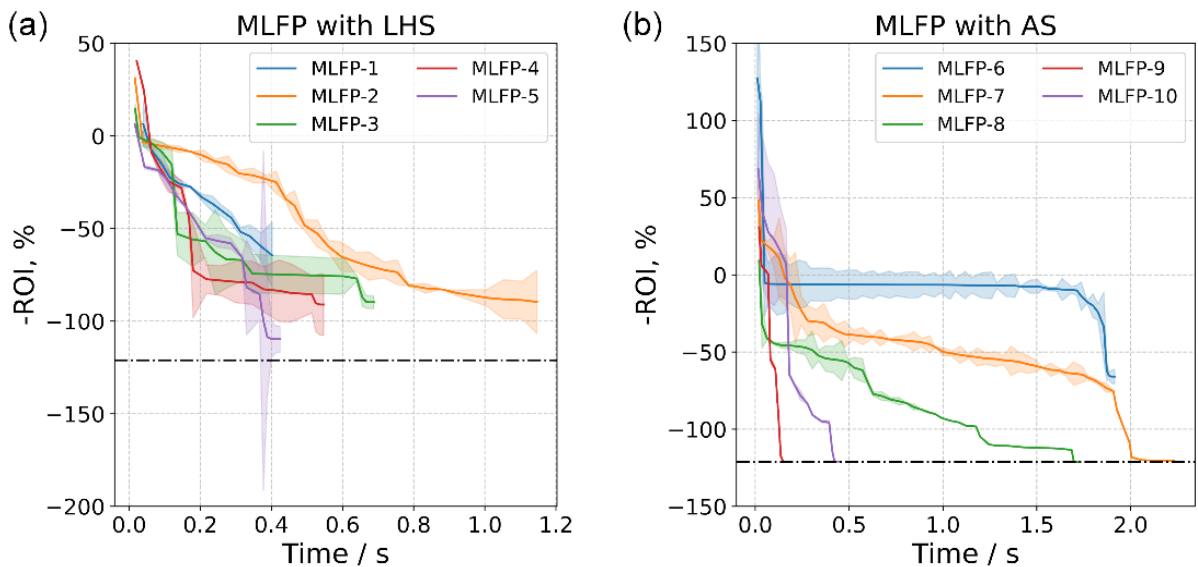

FIGURE 7 Convergence curves for Example 11 with error band of MLFP 1-10.

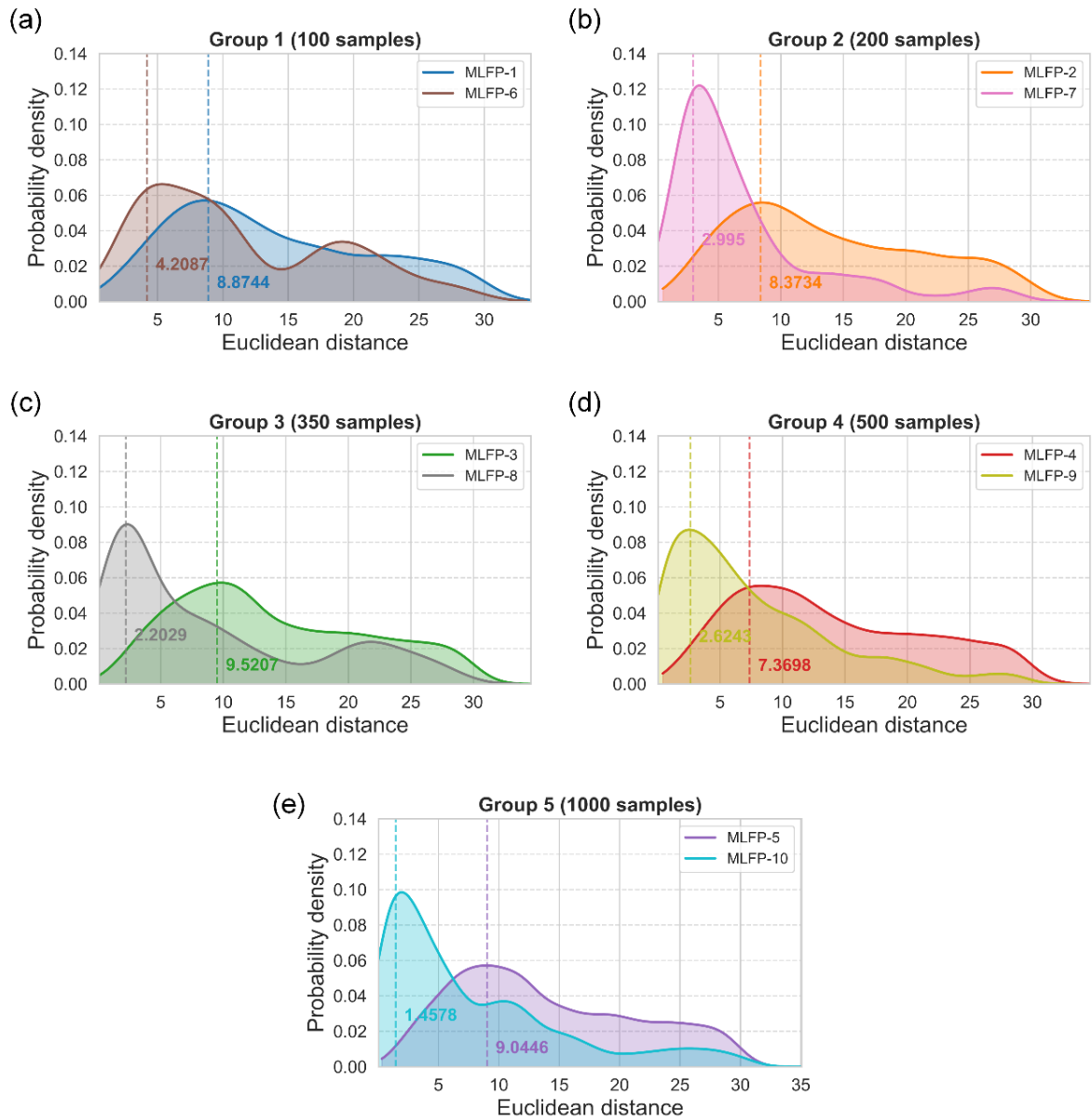

FIGURE 8 Probability density function of Euclidean distance from sampling points to rigorous model-based solution for Example 11.

The comparative results between TRFs 1-4 and MLFPs 1-10 are illustrated in FIGURE 9 where the dashed line represents the optimal solution from the rigorous models solved by IPOPT[81]. Compared to TRFs1-4, MLFPs 7-8 achieve the same optimal ROI of 121.1 but with significantly fewer function evaluations. FIGURE 10 presents the optimization results of TRFs 1–4 and MLFPs 1–10 with 10 runs. As shown in FIGURE 10 (a), MLFPs 6-10 are much more stable than MLFPs 1-5 with the same function evaluations, indicating that MLFPs with the proposed adaptive sample strategy is superior to those with LHS. When the number of samples increases, the stability of MLFPs 6-10 increases, evidenced by the reduced variance in optimized ROI across multiple runs. From FIGURE 9 (a) and FIGURE 10 (b), we can observe that both TRF-4 and MLFP-10 require similar function evaluations (915 vs. 1002) and the rigorous model-based optimal solution. However, MLFP-10 shows greater stability compared to TRF-4 due to the box height (Interquartile Range, IQR) and whisker difference of MLFP-10 are significantly narrower than TRF-4.

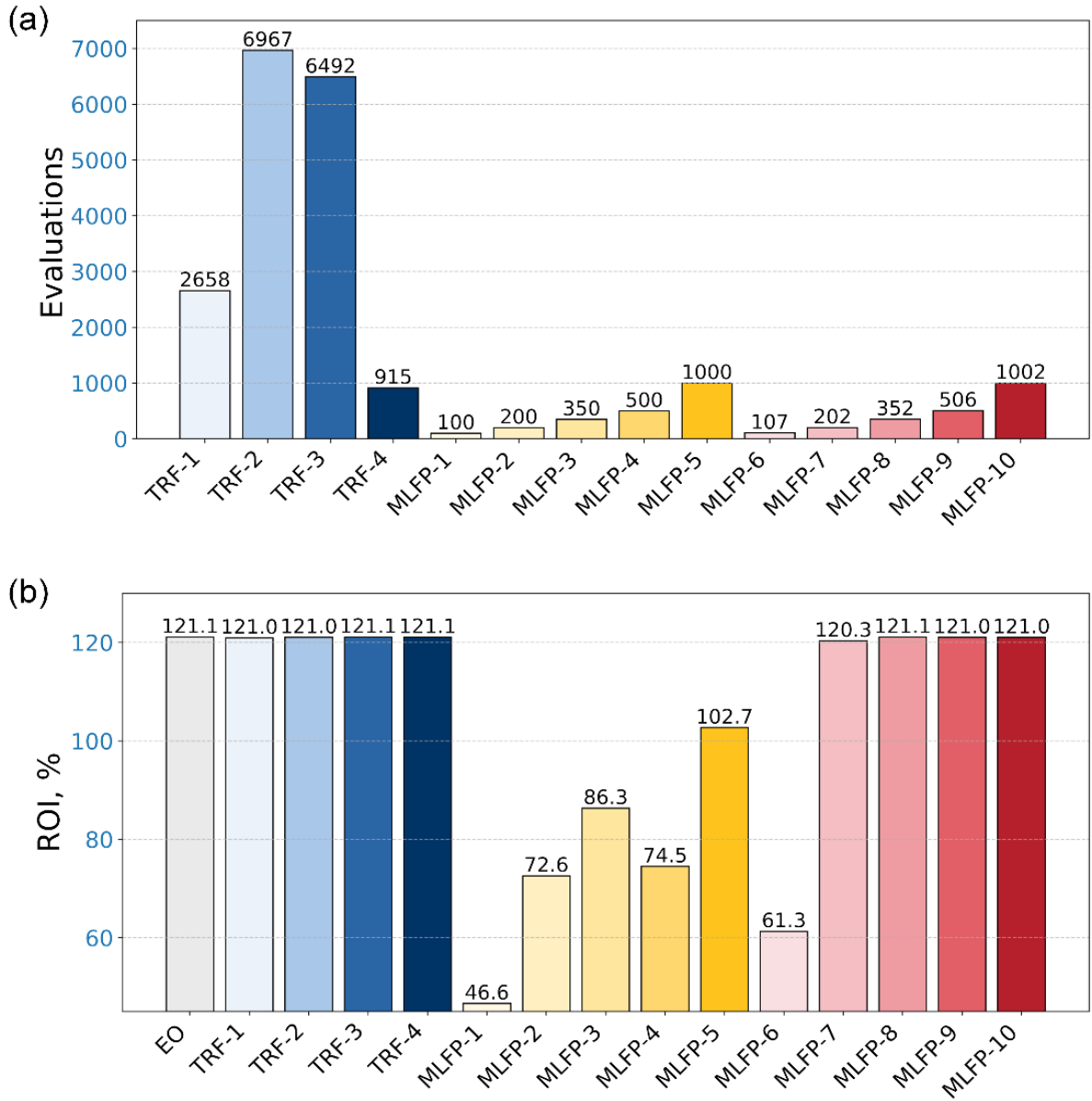

FIGURE 9 ROI results for Example 11 from MLFPs and TRF1-4 on (a) no. of function evaluations of the black-box function and (b) objective function value.

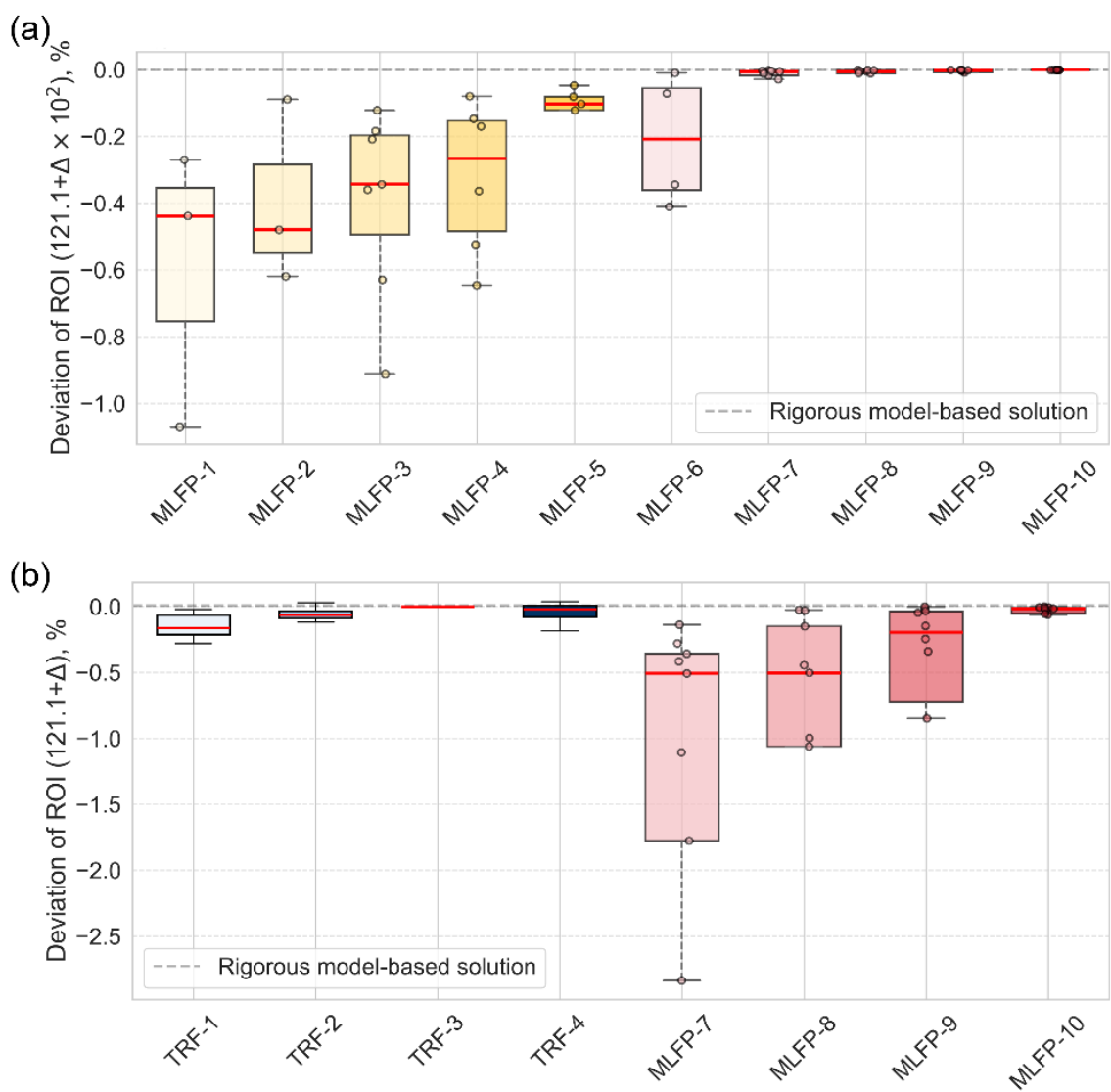

FIGURE 10 Boxplot of best ROI for Example 11 from (a) MLFP1-10, and (b) TRF 1-4 and MLFP 7-10 with 10 runs.

### 4.3 Example 12: Extractive distillation (ED) for toluene and n-heptane separation

Process simulation-based optimization is a typical BBO problem involving expensive function evaluations or unknown functions. In this subsection, we test an example of extractive distillation to separate toluene form n-heptane using solvent phenol provided by Ma et al. (2022). The process flowsheet is provided in FIGURE S7 of the **Supplementary Material**. The process mainly consists of two columns. The first column (C1) is used to separate n-heptane, and the second column (C2) is used for separation of toluene form solvent phenol. The purity requirement for both products is not less than 0.98. All specific constraints and utility prices are referenced from Fan et al.[53]. The optimization problem using the surrogate models mainly achieves the minimum operating cost by adjusting the amount of extractant, reflux ratio, and distillate flow rate, as shown in (**PED**).

$$
\begin{aligned}
&\min_{F, r_1, D_1, r_2, D_2 \in \mathbb{R}} && c_{hu}(Q_2 + Q_4) - c_{cu}(Q_1 + Q_3) \\
&s.t. && (x_n, x_t, Q_1, Q_2, Q_3, Q_4) = s(F, r_1, D_1, r_2, D_2) \\
& && 0.98 \leq x_n \leq 1 \\
& && 0.98 \leq x_t \leq 1
\end{aligned}
\qquad \textbf{(PED)}
$$

where, $c_{hu}$ and $c_{cu}$ is the cost coefficients for cold and hot utilities; $Q_1$ and $Q_3$ are duties of condensers C1 and C2, $Q_2$ and $Q_4$ are duties of reboilers C1 and C2, kW; $F$ is the mass flow of solvent phenol, kg·hr$^{-1}$; $r_1$ and $r_2$ are the reflux ratio of C1 and C2; $D_1$ and $D_2$ are the distillate flow of C1 and C2; $x_n$ and $x_t$ are the mole fraction of n-heptane in C1 distillate and toluene in C2 distillate; $s(\cdot)$ is the surrogate model for extractive distillation process. The variable range and parameter values are obtained from Fan et al. [53].

According to Ma et al.[24,70], the surrogate models constructed using either ReLU neural networks or ALAMO exhibit great relative errors ranging from 0.2% to 25%. Despite employing complex basis functions, these models still fail to fully capture behaviors of the process system due to inherent limitations in their nonlinear representation capabilities. We use the proposed MLFP framework to generate a surrogate model for the entire process where 1014 valid sampling points are generated. The performance of the constructed surrogate models on the test set is provided in FIGURE S8 of the **Supplementary Material**. Each graph exhibits a correlation coefficient of 1.00 ($R^2 \cong 1$) and relative errors are no greater than 0.56%, indicating excellent predictive performance.

We then use MLFP to solve the **PED** problem. The convergence curves of the objective value, and n-heptane and toluene molar fractions are depicted in FIGURE 11. As illustrated in FIGURE 11, the predicted objective value, and n-heptane and tolumen molar fractions from the surrogated models are in good agreement with those validated from Aspen rigorous simulation, indidating the high model prediction accuracy. The best energy cost is 134.60 $ hr$^{-1}$. The achieved purities of n-heptane and toluene products are 0.98 and 0.98

respectively, which satisfy the purity requirements.

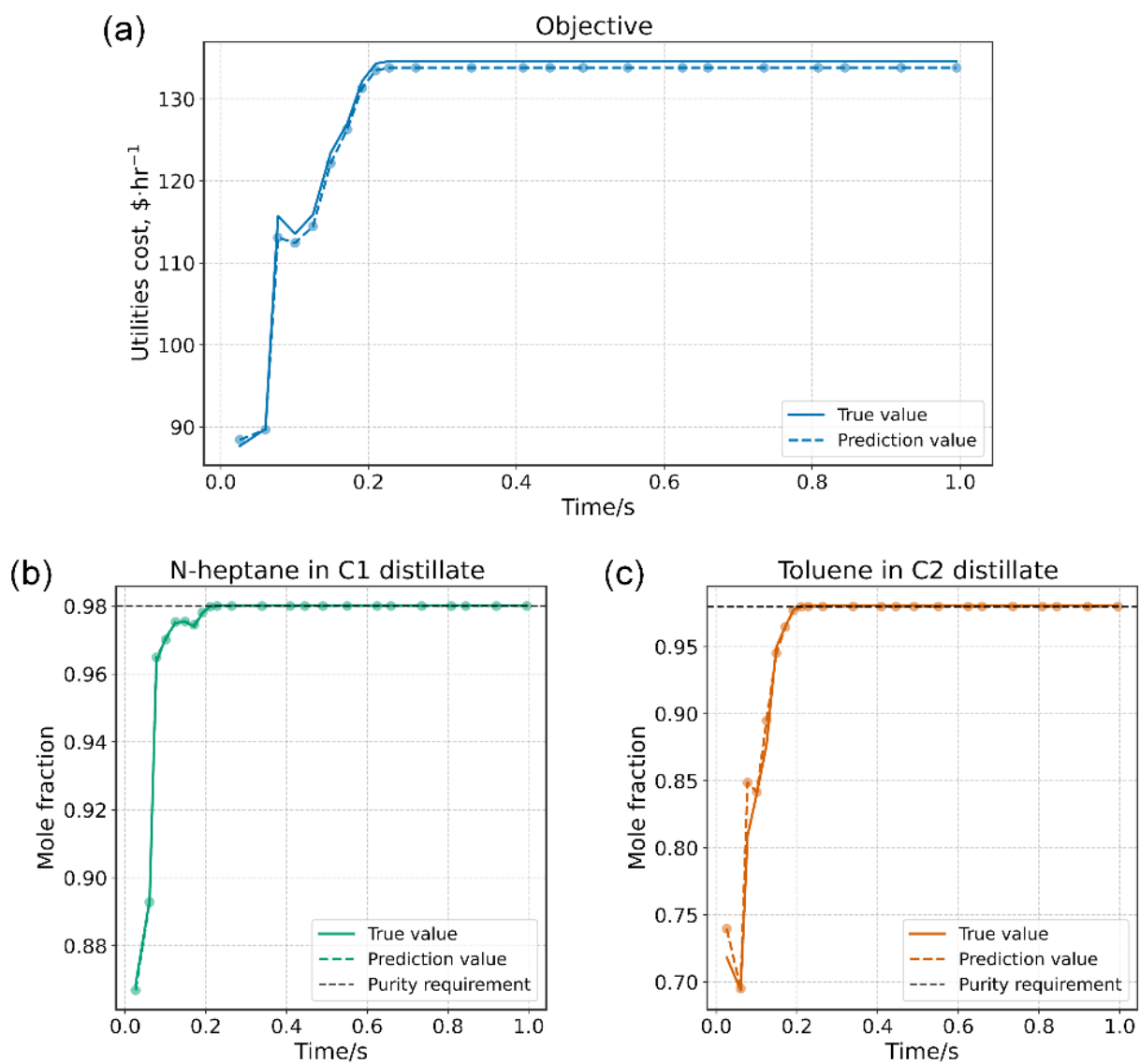

FIGURE 11 Convergence curves for Example 12 where (a) objective values, (b) n-heptane molar fraction in the top distillate of C1 and (c) toluene molar fraction in the top distillate of C2.

To further illustrate the capability of the proposed MLFP framework, we also use genetic algorithm (GA) [82], particle swarm optimization (PSO)[82] to solve problem **PED** where the surrogate model is obtained from the proposed MLFP framework. The reported network structure with the highest accuracy in Ma et al.[70] is also used to develop ReLU NN by using the sampling points generated from LHS, which is then reformulated into an MILP problem and optimized by Gurobi[83]. All optimization results using the surrogate models are validated in Aspen Plus. The validated results are provided in Table 1. The best results from the adaptive sampling Bayesian optimization (ASBO-I)[53] are reported by Tian et al. and also provided in Table 1 for comparison.

TABLE 1 Validated results for Example 12 from different optimization strategies.

| Algorithms | Objective ($\cdot hr^{-1}$) | Purity | | No. of samples/Initial points | Optimization time (sec) | Total time (min) |
|---|---|---|---|---|---|---|
| | | $x_n$ | $x_t$ | | | |
| Surrogate-GA | 154.77 | 0.99 | 0.98 | 2400 | 3.59 | 62.14 |
| Rigorous-GA | 177.07 | 0.98 | 1.00 | - | 2898.55 | 48.31 |
| Surrogate-PSO | 134.73 | 0.98 | 0.98 | 2400 | 3.94 | 62.15 |
| Rigorous-PSO | 151.60 | 0.98 | 0.98 | - | 2919.47 | 48.66 |
| SQP (Aspen) | 231.18 | 0.99 | 0.99 | 1 | Failed | - |
| ASBO-I (best) | 248.60 | valid | valid | 160 | - | - |
| ReLU NN | 154.13 | 0.96 | 0.96 | 2400 | 2.57 | 79.11 |
| MLFP | 134.60 | 0.98 | 0.98 | 2400 | 0.99 | 62.08 |

From TABLE 1, it is shown that SQP (Aspen) fails to converge and halts at the initial point provided. ASBO-I requires the fewest samples, but its best objective is 248.60 $\cdot hr^{-1}$, which is 84.7% higher than that of 134.60 $\cdot hr^{-1}$ obtained from MLFP. Although the reformulated ReLU NN using the network structure reported in Ma et al.[70] obtains the best objective of 154.13 $\cdot hr^{-1}$, which is lower than ASBO-I, but still higher than that of 134.60 $\cdot hr^{-1}$ from MLFP. Moreover, the purity of n-heptane and toluene is 0.96, which does not meet the purity requirements. In contrast, MLFP exhibits the best performance, generating the lowest objective value of 134.60 $\cdot hr^{-1}$ and all purity constraints are met.

## 4.4 Example 13: $CO_2$ capture from biogas (CCB)

Removing $CO_2$ from biogas is an important step in producing biofuels as it can increase the concentration of $CH_4$, thereby enhancing the quality of biogas as an energy carrier. Purified biogas, commonly known as biomethane or renewable natural gas (RNG), can be used for various purposes, such as direct combustion to generate heat or electricity, as vehicle fuel, or injection into natural gas pipelines.

This example is from Xu et al.[80] where n-methyldiethanolamine (MDEA) is used as the solvent. The reactions between $CO_2$ and MDEA take place in the absorber, as given in Eq.(2). The carbonic acid undergoes thermal decomposition in the stripper, as illustrated in

Eq.(3).

$$CO_2 + R_2NCH_3 + H_2O \longrightarrow R_2NCH_3H^+ + HCO_3^-, \qquad (2)$$

$$R_2NCH_3H^+ + HCO_3^- \xrightarrow{\Delta} R_2NCH_3 + CO_2 \uparrow + H_2O. \qquad (3)$$

The desorbed MDEA is circulated in the system as the regenerated solvent. The products $CH_4$ and $CO_2$ will carry a small amount of water vapor and MDEA, so it is necessary to supplement the solvent appropriately. The entire flow diagram can be found in FIGURE S9 of the **Supplementary Material**. We employ Python to externally call Aspen's Application Programming Interface (API) for direct iteration, ensuring that the loss remains consistent with the replenished solvent.

We treat the entire process as a black box. The surrogate model is constructed by using the proposed MLFP framework to predict $CH_4$ molar fraction in the vapor outlet stream of Flash1, $CO_2$ mole fraction in the vapor outlet stream of Flash2, fresh MDEA required, and duties of reboiler based on the Desorber. We then establish the optimization problem, which is denoted as (**PCCB**). The objective is to minimize total annual cost, including utility costs, annualized equipment investment costs, and supplementary MDEA costs. The purified biogas must achieve a minimum $CH_4$ purity of 0.97, while the $CO_2$ content in the regenerated MDEA solution should be less than 0.1.

$$
\begin{aligned}
\min_{\substack{F_{MDEA}^R \in \mathbb{R} \\ N_A, N_D \in \mathbb{R}}} \quad & HrF_{MDEA}^S - c_{cu}Q_c + c_{hu}Q_h + \frac{[c_{MDEA}F_{MDEA}^R + c_s(N_A + N_D)]}{Yr} \\
s.t. \quad & (x_{CH_4}, x_{CO_2}, F_{MDEA}^S, Q_c, Q_h) = s(F_{MDEA}^R, N_A, N_D) \\
& 0.97 \leq x_{CH_4} \leq 1 \\
& 0 \leq x_{CO_2} \leq 0.1 \\
& 2.7 \leq F_{MDEA}^R \leq 5 \\
& 3 \leq N_A \leq 5 \\
& 3 \leq N_D \leq 11 \\
& N_A, N_D \in \mathbb{R}
\end{aligned}
\qquad \text{(PCCB)}
$$

where, $Hr$ is the annual operating time, set to 8000 hours; $Yr$ is the period of depreciation, set to 10 years; $c_{hu}$ and $c_{cu}$ is the energy cost parameters for cold and hot utilities, set to 20 and 80 $\cdot kW^{-1}\cdot yr^{-1}$ [84]; $c_{MDEA}$ is the price of MDEA, set to 13000 $\cdot ton^{-1}$; $c_s$ is the price for each column stage, set to 8508 $ [85,86]; $Q_c$ and $Q_h$ are cooling and heating utilities of the whole system, kW; $F_{MDEA}^R$ is the mass flow of recycled MDEA, ton$\cdot hr^{-1}$; $F_{MDEA}^S$ is the mass flow of supplementary MDEA, ton$\cdot hr^{-1}$; $N_A$ and $N_D$ are the number of stages of absorber and desorber, which is assumed to be continuous; $x_{CH_4}$ and $x_{CO_2}$ are the mass fraction of $CH_4$ and $CO_2$ in system outlets; $s(\cdot)$ is the surrogate model for the entire process simulation.

The performance of the constructed surrogate model on the test set is shown in FIGURE S10 of the **Supplementary Material**, and the relative error on the test set of all predicted outputs is less than 0.1%. The convergence curves for the objective function value, $CH_4$ mass fraction in the purified biogas and $CO_2$ mass fraction in the recycled lean MDEA when solving problem **PCCB** are illustrated in FIGURE 12. As shown in FIGURE 12, the dash line representing the predicting values overlaps with the solid line representing the true values. In other words, the prediction values well match the validated values from Aspen Plus, indicating the high accuracy of the constructed surrogate model. The best annual cost of $2.55\times10^5$ $\cdot yr^{-1}$ is obtained within 0.2 s. The mass fraction of $CH_4$ in the purified biogas is 0.97, which satisfied the purity requirement after validation. Although the number of stages in the absorber and stripper is assumed to be continuous without integrality constraints, the number of stages of both columns is 3 in the final best solution, which is coincidentally integer.

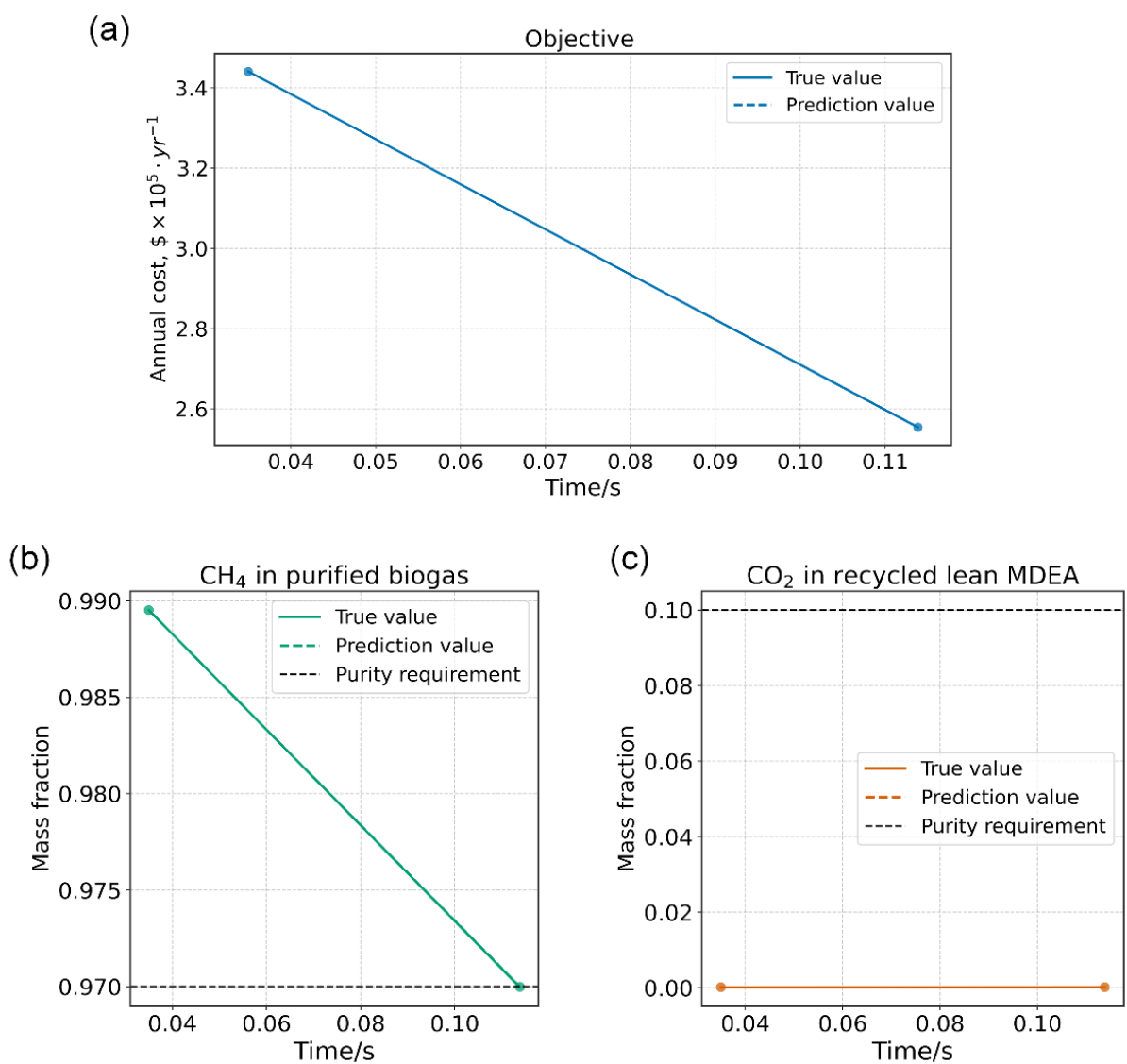

FIGURE 12 Convergence curves for Example 13 (a) objective function value; (b) mass fraction of $CH_4$ in the purified biogas, and (c) mass fraction of $CO_2$ in recycled lean MEDA.

## 5. Conclusion

A machine learning powered feasible path (MIFP) optimisation framework integrated with adaptive sampling was proposed for solving general black-box optimisation (BBO) problems. Through adaptive sampling, the number of samples required to construct the surrogates by machine learning algorithms is significantly reduced, as the proposed strategy could filter out a substantial portion of infeasible samples and automatically identify promising regions that were potentially optimal. By incorporating feasible path method, MLFP managed to avoid the substantial brainy/nonlinear terms, that typically arise from the full-space algebraic formulations of machine learning models, ensuring the rapid convergence.

Computational studies demonstrate that the proposed framework was capable of constructing surrogate models and converging to solutions around the KKT point in Example 1-11, even with a limited number of samples. Compared with the state-of-the-art algorithms, MLFP shows similar performance to TRFs 1-4 for Example 6-10, but requires fewer evaluations; For Example 12, the optimal objective value obtained with MLFP is 84.7% lower than the best result of ASBO-I from its 10 runs. Moreover, multiple independent runs exhibit similar performance, indicating strong robustness of the proposed framework. We will continue to explore applications of MLFP for global optimisation of NLP and MINLP problems in our future research.

**Supporting information**

Additional supporting information can be found online in the Supporting Information section at the end of this article. Supplementary Data for this paper is available at https://data.mendeley.com/preview/rk6jn8wg55/2. The raw data set sampled from Examples are provided in Data A as a .zip file. The Aspen Plus simulations of Example 12 and 13, the EO model of Example 11, and the programs to optimize the EO model of Expamle11 and the ReLU neural network model of Example 12 are provided in Data B as a .zip file. The neural network surrogate models for all Examples are provided in Data C as a .zip file. The other supplementary information is provided in **Supplementary Material** as a .docx file. An available example code of proposed algorithm is provided at https://github.com/ZixuanChang/MLSQP .

**Supplementary Material for**

**Machine Learning Powered Feasible Path Framework with Adaptive**

**Sampling for Black-box Optimization**

**Content**

## S1. Latin hypercube sampling and adaptive sampling

We use an illustrative example to compare the performance of Latin hypercube sampling (LHS) and the proposed adaptive sampling strategy (AS). In this example, the objective function is six-hump camel function (see Eq.(S1)), which is a two-dimensional non-convex function with 2 global minima $-1.0316$ at $(0.0898, -0.7126)$ and $(-0.0898, 0.7126)$. The function is evaluated within $[-1,1]$. To enhance the generality of the example, a constraint $x_2 - x_1 \leq 0$ is added. The objective function and constraint are both considered as black-box functions.

$$f(x_1, x_2) = \left(4 - 2.1 \cdot x_1^2 + \frac{x_1^4}{3}\right)x_1^2 + x_1 \cdot x_2 + 4\left(-1 + x_2^2\right) \cdot x_2^2 \tag{S1}$$

By employing the two-stage training strategy presented in FIGURE 2, AS first collects a small amount of samples (e.g., 50) by LHS, and a support vector machine (SVM) classifier is trained to approximate $x_2 - x_1 \leq 0$. Due to the small number of sample points, the SVM classifier may incorrectly classify feasible regions as infeasible. We then extract the support vectors of the pre-trained SVM and forcefully designate them as feasible points to train a more conservative SVM. The spatial distribution of the samples in the $x_1 - x_2$ plane from LHS and AS is presented in FIGURE S1. As shown in FIGURE S1, LHS sampling points contain 51.3% feasible points, while AS contains 91.23%, which shows that AS can effectively reduce the infeasibility of sampling points.

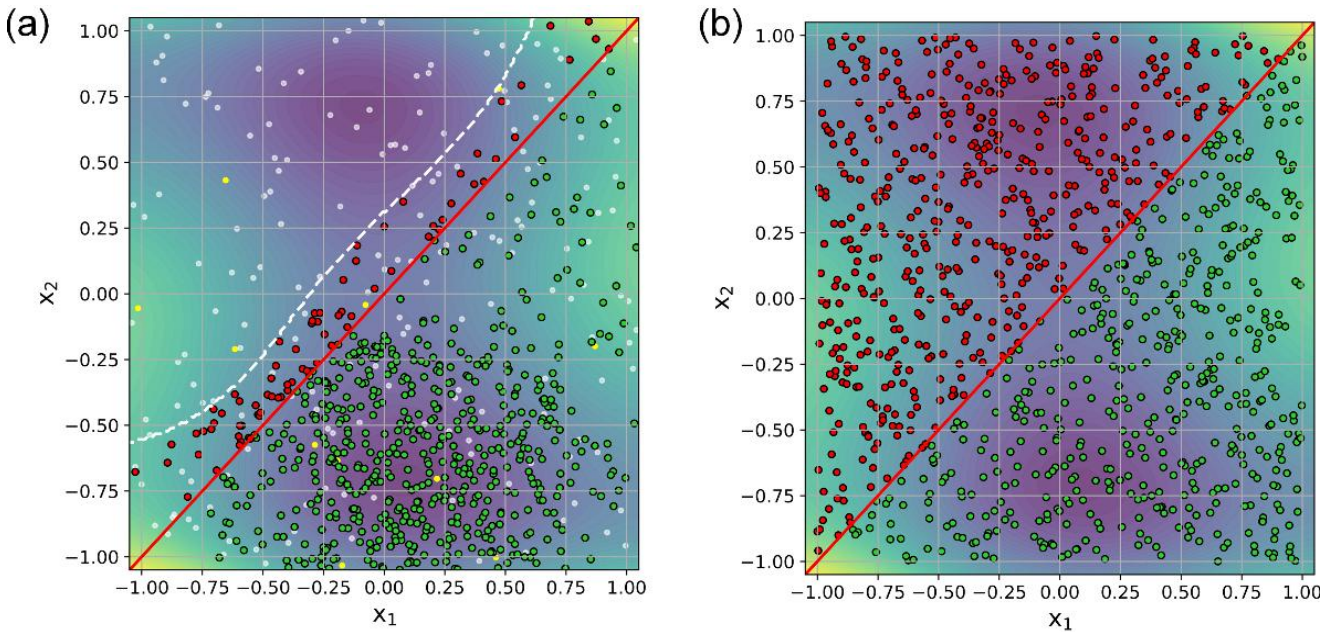

FIGURE S13 Sampling results of (a) AS and (b) LHS where red points reperesent infeasible points, green points denote feasible points, white points reperesent initial LHS points for AS, red solid line represents $x_1 \leq x_2$ plane, white dashed line represents SVM classifier.

## S2. Surrogates with derivatives

Several machine learning algorithms such as decision tree, support vector regression machine, neural networks, etc., can be used to develop surrogate models for black-box functions. In the following, they are introduced, and their derivatives are also explained.

## S2.1. Decision Tree (DT)

Various forms of decision trees have been widely used as surrogate models in optimization[1–3]. In the early stage, the decision boundaries of decision trees are univariate and parallel to the coordinate axes, using the mean as the predictive model. Subsequently, decision trees with linear predictive models[1] are developed. In recent studies, decision trees employ linear hyperplanes as decision boundaries and use arbitrary functions as predictive models to fit sample data[4].

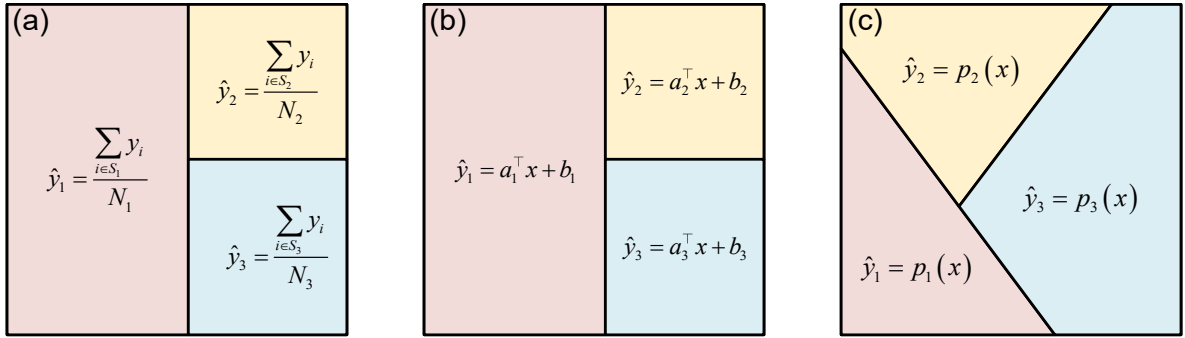

FIGURE S14 Decision tree with (a) single variable decision boundary and mean prediction model, (b) single variable decision boundary and linear prediction model as well as (c) linear hyperplane decision boundary and arbitrary prediction model.

In this work we do not require the specific form of the decision boundary. Instead, we require the predictive model in each leaf is twice differentiable and the second-order derivatives are not all equal to zero. Considering the single-variable decision boundary can be regarded as a special form of linear decision boundary, the generalized algebraic formulation can be stated as follows,

$$\hat{y}^{DT} = p_i(x) \qquad \forall x \in \mathcal{P}_i, i \in \mathcal{L}, \qquad (S2)$$

$$\mathcal{P}_i = \left\{ x \in \mathbb{R}^n : a_j^\top x \le b_j, \forall j \in L(L_i); x \in \mathbb{R}^n : a_j^\top x > b_j, \forall j \in R(L_i) \right\} \qquad (S3)$$

where $\mathcal{L}$ is the set of leaves; $\hat{y}^{DT}$ is the predictive value of decision tree, $x$ is the input feature; $L_i$ represents the leaf node $i$, $j$ represents the branch node; $p_i(\cdot)$ is the the predictive model at $L_i$. For a binary tree, the path from the root to the leaves exists and is unique. For branch node $j$ in the path from the root to leaf $L_i$, as is shown in FIGURE S15 (a), if $x$ satisfies $a_j^\top x \leq b_j$, then the node $j$ belongs to $L(L_i)$, otherwise belongs to $R(L_i)$. $\mathcal{P}_i$ represents the polyhedron as is stated in Eq.(S3).

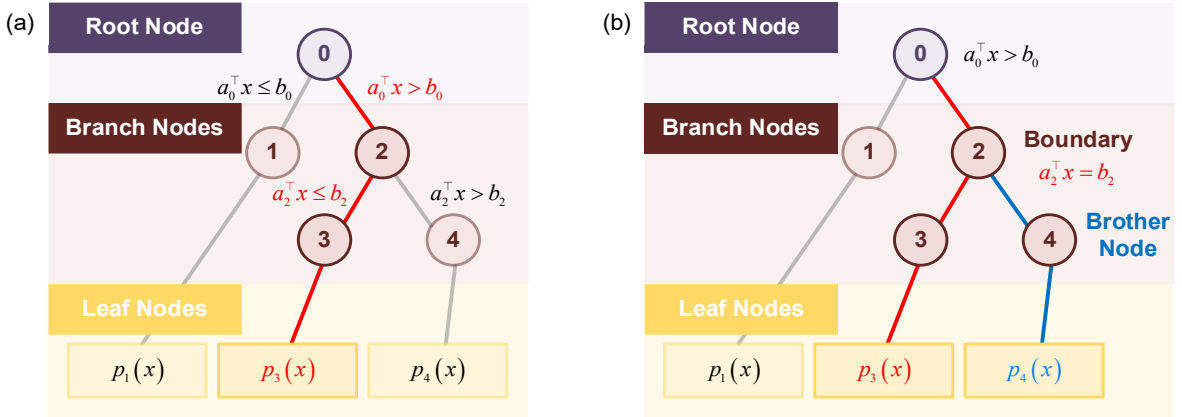

FIGURE S15 (a) The path from the root to the leaves and (b) the boundary prediction functions' recognition.

For those points that do not exactly satisfy the equality boundary, the first- and second-order derivatives are as follows:

$$\nabla_x \hat{y}_{DT} = \nabla_x p_i(x) \quad \forall x \in \mathcal{P}_i, i \in \mathcal{L} \tag{S4}$$

$$\nabla_x^2 \hat{y}_{DT} = \nabla_x^2 p_i(x) \quad \forall x \in \mathcal{P}_i, i \in \mathcal{L} \tag{S5}$$

For those points that exactly satisfy the equality boundary, they are non-differentiable. Therefore, the minimum norm subgradient is used in place of the first-order derivative:

$$\partial \hat{y}_{DT} = \min_{\omega_i} \left\| \sum_i \omega_i \nabla_x p_i(x) \right\|_2 \quad x \in \bigcup_{j \in L(L_i) \cup R(L_i)} \left\{ a_j^\top x = b_j \right\} \tag{S6}$$

Then, the finite difference method is employed to estimate the second-order derivative and construct an approximate Hessian matrix. The algorithm for boundary nodes' identification and derivatives' calculation of decision tree is shown in ALGORITHM S1.

ALGORITHM S1 Algorithm for boundary nodes' identification and derivatives' calculation of decision tree

**Data:** Decision Tree $\mathcal{T}$

**Input:** $x$, $\mathcal{N}_{eq} \leftarrow \emptyset$, $\mathcal{N}_{path} \leftarrow \emptyset$, $\mathcal{N}_{leaf} \leftarrow \emptyset$

**Output:** $\nabla_x y, \nabla_x^2 y$

```
 1  begin
 2      function cal_path (x, T_ini): {
 3      while True do
 4          N_path.append(T_i)
 5          if a_i^T x − b_i ≤ 0 then
 6              T_{i+1} ← left child node of T_i
 7              if a_i^T x − b_i = 0 then
 8                  N_eq.append(T_i)
 9              end if
10          else
11              T_{i+1} ← right child node of T_i
12          end if
13          if T_{i+1} is leaf node then
14              N_leaf.append(T_{i+1})
15              break
16          end if
17      end while
18      }
19      cal_path(x, T_0)
20      if N_eq ∩ N_path = ∅ then
21          return ∇_x y, ∇_x^2 y calculated using Eq.(5) and Eq.(6)
22      else
23          while N_eq ∩ N_path != ∅ do
24              N_co = N_eq ∩ N_path
25              for n in N_co do
26                  if a_i^T x − b_i ≤ 0 then
27                      T_ini ← right child node of T_n
28                  else
29                      T_ini ← left child node of T_n
30                  end if
31                  cal_path(x, T_ini)
32                  N_eq.delete(T_n)
33              end for
34          end while
35      end if
36      return ∇_x y calculated using Eq.(7) and ∇_x^2 y using finite difference
37  end
```

## S2.2 Support Vector Regression Machine (SVR)

SVR is the application of Support Vector Machine (SVM) in regression problems. SVR seeks to identify a function that best fits the data while maintaining a certain tolerance for error, which is achieved through the introduction of an ε-insensitive loss function[5]. This model is especially effective for handling high-dimensional data and can construct non-linear decision boundaries in the feature space when kernel methods are utilized.

The algebraic form of a kernel SVR is shown in Eq.(S7):

$$\hat{y}_{SVR} = \sum_i \left( \alpha_i - \alpha_i^* \right) \mathcal{K}\left( x, x_i \right) + \beta \tag{S7}$$

where $\hat{y}_{SVR}$ is the predictive value of SVR, $\alpha_i$ and $\alpha_i^*$ is the Lagrange multiplier, $\beta$ is the bias term, $x_i$ is the sample $i$ in the training set, $x$ is the input feature and $\mathcal{K}(\cdot)$ is the kernel function.

And the derivatives are as shown in Eqs.(S8) and (S9).

$$\nabla_x \hat{y}_{SVR} = \sum_i \left( \alpha_i - \alpha_i^* \right) \nabla_x \mathcal{K}\left( x, x_i \right) \tag{S8}$$

$$\nabla_x^2 \hat{y}_{SVR} = \sum_i \left( \alpha_i - \alpha_i^* \right) \nabla_x^2 \mathcal{K}\left( x, x_i \right) \tag{S9}$$

## S2.3 Neural Networks (NN)

Neural networks serve as the cornerstone of deep learning technology, and their applications within chemical engineering are expanding rapidly. This growth is primarily attributed to neural networks' robust nonlinear modeling capabilities and their proficiency in processing intricate data.

By leveraging these strengths, neural networks enable advanced process optimization, predictive maintenance, and innovative product development in the chemical industry, thereby transforming traditional practices with more efficient and intelligent solutions.

The most basic NN model is the multi-layer perceptron (MLP). The algebraic form of MLP with $L$ hidden layers can be expressed by the following composite function:

$$\hat{y}_{MLP} = f^{(L+1)} \circ \phi^{(L)} \circ f^{(L)} \circ \cdots \circ \phi^{(1)} \circ f^{(1)}(x) \tag{S10}$$

Where, $\hat{y}_{MLP}$ is the predictive value of MLP, $f^{(n)}(z) = W^{(n)}z + b^{(n)}$, $\phi^{(n)}(\cdot)$ is the activation function of the n-th layer. $W^{(n)}$ and $b^{(n)}$ represent the weights and biases from layer $n-1$ to layer $n$, respectively.

According to the chain differentiation rule, we have:

$$\frac{\partial \hat{y}_{MLP}}{\partial x} = \frac{\partial f^{(L+1)}}{\partial \phi^{(L)}} \cdot \frac{\partial \phi^{(L)}}{\partial f^{(L)}} \cdots \frac{\partial \phi^{(1)}}{\partial f^{(1)}} \cdot \frac{\partial f^{(1)}}{\partial x} \tag{S11}$$

Specifically, for every layer we have:

$$\frac{\partial f^{(n+1)}}{\partial \phi^{(n)}} = W^{(n+1)} \tag{S12}$$

$$\frac{\partial \phi^{(n)}}{\partial f^{(n)}} = diag\left(\phi'^{(n)}\left(f^{(n)} \circ \phi^{(n-1)} \circ f^{(n-1)} \circ \cdots \circ \phi^{(1)} \circ f^{(1)}(x)\right)\right) \tag{S13}$$

Therefore, the first-order derivative can be expressed as follows:

$$\frac{\partial \hat{y}_{MLP}}{\partial x} = W^{(L+1)} \prod_{n=L}^{1} diag\left(\phi'^{(n)}\left(f^{(n)} \circ \phi^{(n-1)} \circ f^{(n-1)} \circ \cdots \circ \phi^{(1)} \circ f^{(1)}(x)\right)\right) W^{(n)} \tag{S14}$$

To simplify the expression of second-order derivative of MLP, we assume that the predicted output value $\hat{y}_{MLP}$ is a scalar, and the conclusions obtained can be easily extended to high-dimensional output neural networks.

From the output layer to the hidden layer $L$, we have:

$$\frac{\partial^2 \hat{y}_{MLP}}{\partial \mathbf{x} \partial \mathbf{x}^\top} = W^{(L+1)} \cdot \frac{\partial^2 \phi^{(L)}}{\partial \mathbf{x} \partial \mathbf{x}^\top} \tag{S15}$$

For hidden layer $l$ we have:

$$\frac{\partial^2 \phi^{(l)}}{\partial \mathbf{x} \partial \mathbf{x}^\top} = \phi'''^{(l)}\left(f^{(l)}\right) \cdot \left(\frac{\partial f^{(l)}}{\partial \mathbf{x}}\right) \cdot \left(\frac{\partial f^{(l)}}{\partial \mathbf{x}}\right)^\top + \phi'^{(l)}\left(f^{(l)}\right) \cdot \frac{\partial^2 f^{(l)}}{\partial \mathbf{x} \partial \mathbf{x}^\top} \tag{S16}$$

where,

$$\frac{\partial f^{(l)}}{\partial \mathbf{x}} = \prod_{m=1}^{l-1}\left(diag\left(\phi'^{(m)}\left(f^{(m)}\right)\right)W^{(m)}\right) \tag{S17}$$

$$\frac{\partial^2 f^{(l)}}{\partial \mathbf{x} \partial \mathbf{x}^\top} = \sum_{m=1}^{l-1}\left[\left(\prod_{r=m+1}^{l-1}\left(diag\left(\phi'^{(r)}\left(f^{(r)}\right)\right)W^{(r)}\right)\right) \cdot \frac{\partial^2 \phi^{(m)}}{\partial \mathbf{x} \partial \mathbf{x}^\top} \cdot \left(\prod_{k=1}^{m-1}\left(diag\left(\phi'^{(k)}\left(f^{(k)}\right)\right)W^{(k)}\right)\right)\right] \tag{S18}$$

By applying Eqs.(S15)-(S18) and performing layer-by-layer recursion, we can obtain the second-order derivative of the MLP. Obviously, the algebraic symbolic differentiation approach yields expressions that are overly complex and challenging to implement. Therefore, we provide an alternative method: reverse-mode automatic differentiation (AD)[6,7]. Different from central differencing and other numerical differentiation techniques, automatic differentiation provides exact derivatives while significantly reducing computational complexity.

Assuming we have a neural network as shown in FIGURE S16 (a), and the forward process can be transformed into a computational graph as is shown in FIGURE S16 (b). Here, $x_m$ is the input feature of NN, $y$ is the prediction value. $v_i$ and $v_k$ represent the linear transformation in NN, for example $v_1 = w_1 \cdot x_1$, where $w$ represent the weights. $p_j$ represents the activation function, for instance $p_1 = \phi(v_1 + v_4 + b)$, where $b$ represents the bias. By transforming NN into computational graph, the forward trace can be distinctly presented and the backward trace is used for calculating first-order derivative. Similarly, this process can also be transformed into a computational graph used for calculating second-order derivatives. In this work, we use PyTorch[8] to easily achieve automatic differentiation without the tedious derivation for algebraic expressions.

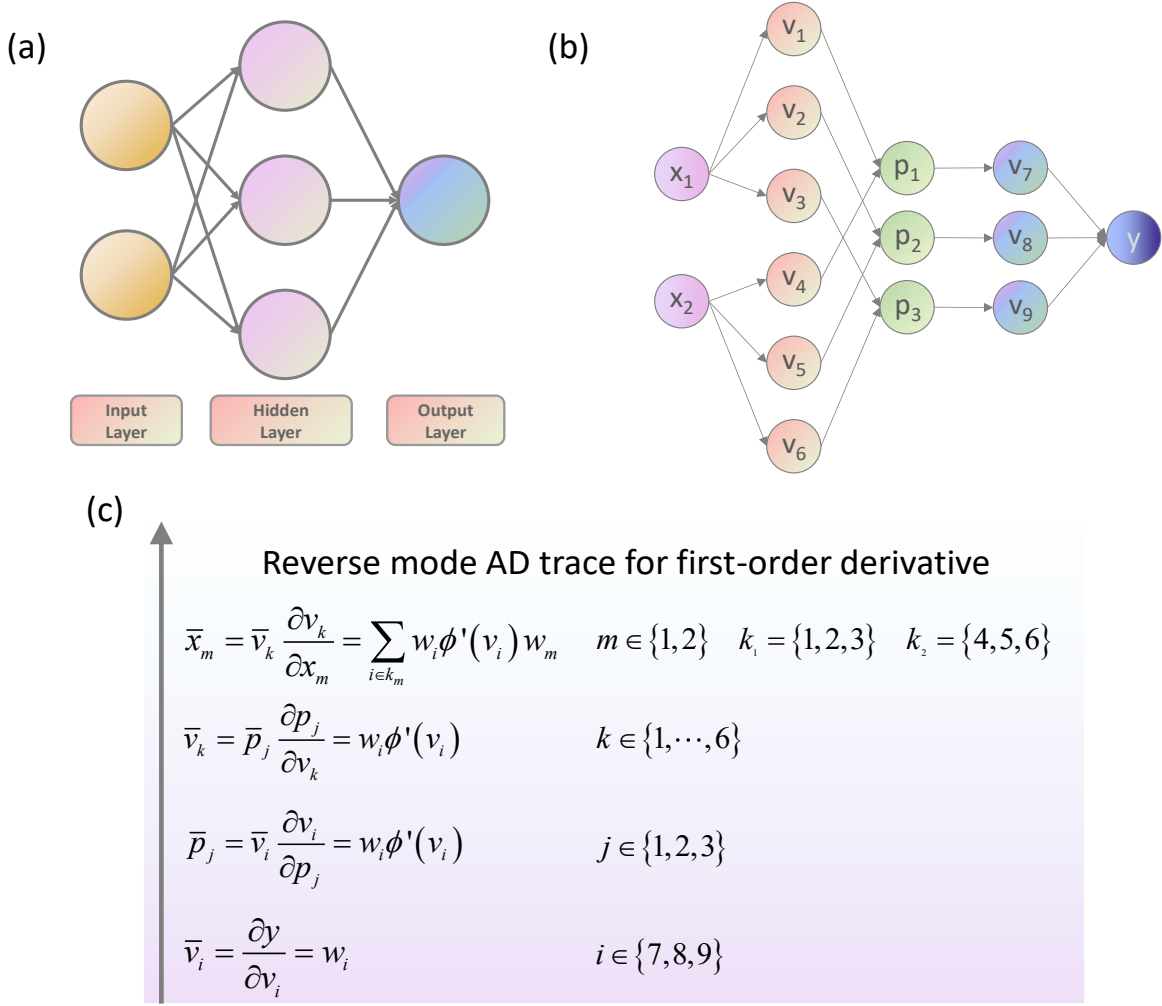

FIGURE S16 An illustrative example of NN with its (a) structure, (b) computational graph and (c) reverse mode AD trace for first-order derivative.

## S 2.4 Ensemble Learning

For some special tasks, a single learner may not perform well. Therefore, the multiple base learners need to be integrated. The weighting method is the most commonly used ensemble approach in regression tasks:

$$\hat{y}_E = \sum_i \omega_i \hat{y}_i \qquad (S19)$$

Where, $\hat{y}_E$ is the prediction value of ensemble learner, $\hat{y}_i$ is the prediction value of base learners and $\omega_i$ is the weights. The derivatives of integrated learner are also the weighted derivatives of the base learner.

The methods outlined above enable the provision of necessary derivative information while preserving the powerful fitting and rapid inference capabilities of machine learning models, thereby facilitating efficient optimization.

**S2.5 Gaussian process regression (GPR)**

GPR, a non-parametric probabilistic regression method, is commonly employed in BO. Given an observation point $x_*$, if we consider the simple special case where the observations are noise free, GP predicts a random variable whose probability distribution follows a normal distribution $\mathcal{N}$. According to definition, the training label $Y$ and the predicted variable $y_*$ follow the following distribution:

$$Y \sim \mathcal{N}\left(0, K\left(X, X\right) + \sigma_n^2 I\right) \tag{S20}$$

$$y_* \sim \mathcal{N}\left[0, \mathcal{K}\left(x_*, x_*\right)\right] \tag{S21}$$

where $\sigma_n^2 \sigma_n^2$ represents the variance of noise, $I$ is the identity matrix, $K(\cdot)$ is the covariance matrix, $\mathcal{K}(\cdot)$ is the kernel function. The joint probability distribution of $Y$ and $y_*$ is a joint normal distribution in the following form:

$$\begin{bmatrix} Y \\ y_* \end{bmatrix} \sim \mathcal{N}\left(0, \begin{bmatrix} K\left(X, X\right) + \sigma_n^2 I & K\left(X, x_*\right) \\ K\left(x_*, X\right) & \mathcal{K}\left(x_*, x_*\right) \end{bmatrix}\right). \tag{S22}$$

GPR uses the expectation of random variables for prediction and measures uncertainty using variance:

$$\hat{y}_* = \mu\left(x_*\right) = K\left(x_*, X\right)^\top \left[K\left(X, X\right) + \sigma_n^2 I\right]^{-1} Y, \tag{S23}$$

$$\sigma^2\left(x_*\right) = \mathcal{K}\left(x_*, x_*\right) - K\left(x_*, X\right)^\top \left[K\left(X, X\right) + \sigma_n^2 I\right]^{-1} K\left(x_*, X\right). \tag{S24}$$

In the BO framework, the problem of maximizing the acquisition function should be resolved before each sampling:

$$x_* = \arg\max_{x \in \mathbb{R}^n} \quad Aq\big(\mu(x), \sigma^2(x)\big)$$

$$\text{s.t.} \qquad h(x) = 0 \qquad \qquad \text{(MAF)}$$

$$g(x) \leq 0$$

where $Aq(\cdot)$ is the acquisition function, such as the upper confidence bound (UCB) and expected improvement (EI); $h(\cdot)$ is the equality constraints and $g(\cdot)$ is the inequality constraints. It should be noted that if the constraint function is also a black-box function, we can use any of the models mentioned above for surrogate[9] and then optimize using the proposed algorithm. Therefore, the derivatives of $\mu(x)$ and $\sigma^2(x)$ are required, and derived as follows:

$$\nabla_x \mu(x) = \frac{\partial K(x,X)^\top}{\partial x}\Big[K(X,X) + \sigma_n^2 I\Big]^{-1} Y \tag{S25}$$

$$\nabla_x \sigma^2(x) = \frac{\partial \mathcal{K}(x,x)}{\partial x} - 2 \cdot \frac{\partial K(x,X)^\top}{\partial x}\Big[K(X,X) + \sigma_n^2 I\Big]^{-1} K(x,X) \tag{S26}$$

If $x$ is an $n_x$-dimensional vector, and there is a total of $n_t$ pieces of training data, then:

$$\frac{\partial K(x,X)^\top}{\partial x} = \begin{bmatrix} \dfrac{\partial \mathcal{K}(x,X_1)}{\partial x_1} & \dfrac{\partial \mathcal{K}(x,X_1)}{\partial x_2} & \cdots & \dfrac{\partial \mathcal{K}(x,X_1)}{\partial x_{n_x}} \\[2mm] \dfrac{\partial \mathcal{K}(x,X_2)}{\partial x_1} & \dfrac{\partial \mathcal{K}(x,X_2)}{\partial x_2} & \cdots & \dfrac{\partial \mathcal{K}(x,X_2)}{\partial x_{n_x}} \\[2mm] \vdots & \vdots & \ddots & \vdots \\[2mm] \dfrac{\partial \mathcal{K}(x,X_{n_t})}{\partial x_1} & \dfrac{\partial \mathcal{K}(x,X_{n_t})}{\partial x_2} & \cdots & \dfrac{\partial \mathcal{K}(x,X_{n_t})}{\partial x_{n_x}} \end{bmatrix}. \tag{S27}$$

The second-order derivative takes the following form:

$$\Big[\nabla_x^2 \mu(x)\Big]_{i,j} = \frac{\partial^2 \mu(x)}{\partial x_i \partial x_j} = \frac{\partial^2 K(x,X)^\top}{\partial x_i \partial x_j}\Big[K(X,X) + \sigma_n^2 I\Big]^{-1} Y \tag{S28}$$

$$\Big[\nabla_x^2 \sigma^2(x)\Big]_{i,j} = \frac{\partial^2 \mathcal{K}(x,x)}{\partial x_i \partial x_j} - 2\frac{\partial K(x,X)^\top}{\partial x_i} A \frac{\partial K(x,X)^\top}{\partial x_j} - 2\frac{\partial^2 K(x,X)^\top}{\partial x_i \partial x_j} A K(x,X) \tag{S29}$$

where $A = [K(X,X) + \sigma_n^2 I]^{-1} Y$.

## S3. SQP-driven feasible path algorithm

Once the derivative information is obtained, the SQP algorithm can be easily applied to surrogate assisted BBO problems. The Lagrangian function of (**P2**) is shown as Eq.(S30):

$$L\left(\mathbf{x},\mathbf{s}(\mathbf{x}),\boldsymbol{\lambda},\boldsymbol{\mu}\right)=f\left(\mathbf{x},\mathbf{s}(\mathbf{x})\right)+\boldsymbol{\lambda}^{\top}\mathbf{h}\left(\mathbf{x},\mathbf{s}(\mathbf{x})\right)-\boldsymbol{\mu}^{\top}\mathbf{g}\left(\mathbf{x},\mathbf{s}(\mathbf{x})\right) \tag{S30}$$

where $\boldsymbol{\lambda}$ and $\boldsymbol{\mu}$ are the Lagrange multipliers for equality and inequality constraints, respectively.

In this work, we take $\mathbf{H}^{k} := \nabla_{\mathbf{x}}^{2} L(\mathbf{x}^{k}, \mathbf{s}(\mathbf{x}^{k}), \boldsymbol{\lambda}^{k}, \boldsymbol{\mu}^{k})$ as the Hessian matrix of the Lagrangian function at $\mathbf{x}^{k}$, and $\mathbf{B}^{k} := \mathbf{H}^{k} + \mathbf{E}^{k}$ as the modified Hessian matrix. $\mathbf{E}^{k}$ denotes the minimum modification under the Frobenius norm. In order to obtain $\mathbf{E}^{k}$, first perform spectral decomposition on $\mathbf{H}^{k}$, so that $\mathbf{H}^{k} = \mathbf{Q}\boldsymbol{\Lambda}\mathbf{Q}^{\top}$, where $\mathbf{Q}$ is a matrix composed of eigenvectors of $\mathbf{H}^{k}$ and $\boldsymbol{\Lambda}$ is a diagonal matrix composed of eigenvalues of $\mathbf{H}^{k}$. In this way, we take $\delta$ as a constant approaching 0, and by modifying the eigenvalues of $\mathbf{H}^{k}$, we adjust all eigenvalues who are smaller than $\delta$ to $\delta$, thus ensuring the positive definiteness of $\mathbf{B}^{k}$. Therefore, we can define $\mathbf{E}^{k} := \mathbf{Q}\mathrm{diag}(\tau_{i})\mathbf{Q}^{\top}$, where

$$\tau_{i} = \begin{cases} 0 & \sigma_{i} > \delta \\ \delta - \sigma_{i} & \sigma_{i} < \delta \end{cases},$$

and $\sigma_{i}$ represents the eigenvalues of $\mathbf{H}^{k}$.

For surrogates with second-order derivative constant at 0, we apply BFGS to approximate $\mathbf{H}^{k}$. The BFGS update formula[10,11] is as follows:

$$\mathbf{B}^{k+1} = \begin{cases} \left(\mathbf{I} - \dfrac{\Delta\mathbf{x} \cdot \Delta L^{\top}}{\Delta L^{\top} \cdot \Delta\mathbf{x}}\right)\mathbf{B}^{k}\left(\mathbf{I} - \dfrac{\Delta L \cdot \Delta\mathbf{x}^{\top}}{\Delta L^{\top} \cdot \Delta\mathbf{x}}\right) + \dfrac{\Delta\mathbf{x} \cdot \Delta\mathbf{x}^{\top}}{\Delta L g^{\top} \cdot \Delta\mathbf{x}} & \text{if } \Delta L^{\top} \cdot \Delta\mathbf{x} > \epsilon \|\nabla_{k} L\| \Delta\mathbf{x}^{\top} \cdot \Delta\mathbf{x}, \epsilon = 10^{-6} \\ \mathbf{B}^{k} & \text{otherwise} \end{cases} \tag{S31}$$

where, $\Delta L := \nabla_{\mathbf{x}} L(\mathbf{x}^{k+1}, \mathbf{s}(\mathbf{x}^{k+1}), \boldsymbol{\lambda}^{k+1}, \boldsymbol{\mu}^{k+1}) - \nabla_{\mathbf{x}} L(\mathbf{x}^{k}, \mathbf{s}(\mathbf{x}^{k}), \boldsymbol{\lambda}^{k}, \boldsymbol{\mu}^{k})$, $\Delta\mathbf{x} := \mathbf{x}^{k+1} - \mathbf{x}^{k}$, and $\mathbf{B}^{0} = \mathbf{I}$.

In this work, we use SQP to solve a surrogate-assisted BBO problem, who expands the original problem into a convex quadratic programming problem through Taylor expansion at each iteration point and conduct line search to merit function along the direction obtained from problem (**QP**). SQP is a type of Newton's method that can effectively handle both small and large-scale nonlinear constrained optimization problems. It can be integrated into line search or trust region frameworks. In this work, we use the line search SQP. Perform second-order Taylor expansion on the Lagrangian function of (**P2**) at $\mathbf{x}^k$ and first-order Taylor expansion on the constraints to obtain the following quadratic programming problem (**QP**):

$$\begin{aligned}
\min_{\mathbf{x} \in \mathbb{R}^n} \quad & \frac{1}{2}\mathbf{d}^{k\top}\mathbf{B}^k\mathbf{d}^k + \nabla f\left(\mathbf{x}^k, \mathbf{s}(\mathbf{x}^k)\right)^{\top}\mathbf{d}^k \\
\text{s.t.} \quad & \nabla \mathbf{h}\left(\mathbf{x}^k, \mathbf{s}(\mathbf{x}^k)\right)^{\top}\mathbf{d}^k + \mathbf{h}\left(\mathbf{x}^k, \mathbf{s}(\mathbf{x}^k)\right) = 0 \\
& \nabla \mathbf{g}\left(\mathbf{x}^k, \mathbf{s}(\mathbf{x}^k)\right)^{\top}\mathbf{d}^k + \mathbf{g}\left(\mathbf{x}^k, \mathbf{s}(\mathbf{x}^k)\right) \leq 0
\end{aligned} \quad \textbf{(QP)}$$

where $\mathbf{d}^k := (\mathbf{x} - \mathbf{x}^k)$ is the search direction of (**P2**), and $\nabla f$, $\nabla \mathbf{h}$ and $\nabla \mathbf{g}$ are the gradients of $f$, $\mathbf{h}$ and $\mathbf{g}$ respectively. The derivative with respect to $x_i$ is $\nabla_{x_i} F(\mathbf{x}, \mathbf{s}(\mathbf{x})) = \frac{\partial F}{\partial x_i} + \sum_{j=1}^{m} \frac{\partial F}{\partial s_j(\mathbf{x})} \frac{\partial s_j(\mathbf{x})}{\partial x_i}$, where $F$ can be any one of $f$, $h$ and $g$.

After solving the problem (**QP**), we obtain the search direction $\mathbf{d}^k$, multipliers $\boldsymbol{\lambda}_{qp}^k$ and $\boldsymbol{\mu}_{qp}^k$ of problem (**QP**) at $\mathbf{x}^k$. Then, line search will be executed on the merit function and the step size $\alpha$ will be determined through backtracking method. In next iteration, $\mathbf{x}^{k+1} = \mathbf{x}^k + \alpha\mathbf{d}^k$ and $\mathbf{H}^{k+1} := \nabla_\mathbf{x}^2 L\left(\mathbf{x}^{k+1}, \mathbf{s}(\mathbf{x}^{k+1}), \boldsymbol{\lambda}_{qp}^k, \boldsymbol{\mu}_{qp}^k\right)$, according to the relationship as the equivalence between SQP and Newton's method [10]. The (**QP**) subproblem is solved by utilizing the qpsolver [12] to call the open-source quadratic programming solver proxqp[13,14].

In this work, the $\ell_1$ merit function, defined as Eq.(S32), is used to evaluate convergence criteria:

$$\phi_1\left(\mathbf{x}^k; \boldsymbol{\lambda}_{qp}^k, \boldsymbol{\mu}_{qp}^k\right) = f\left(\mathbf{x}^k\right) + \left(\boldsymbol{\rho}^k\right)^{\top}\left|\mathbf{h}\left(\mathbf{x}^k\right)\right| + \left(\mathbf{v}^k\right)^{\top}\mathbf{g}\left(\mathbf{x}^k\right)^{-} \tag{S32}$$

where, $\mathbf{g}(\mathbf{x}^k)^{-} := \max\left(0, -\mathbf{g}(\mathbf{x}^k)\right)$, and the penalty parameters are defined as follows[15]:

$$\boldsymbol{\rho}^k = \max\left(\left|\boldsymbol{\lambda}_{qp}^k\right|, \frac{\boldsymbol{\rho}^{k-1} + \left|\boldsymbol{\lambda}_{qp}^k\right|}{2}\right) \tag{S33}$$

$$\mathbf{v}^k = \max\left(\left|\boldsymbol{\mu}_{qp}^k\right|, \frac{\mathbf{v}^{k-1} + \left|\boldsymbol{\mu}_{qp}^k\right|}{2}\right) \tag{S34}$$

Given the $\ell_1$ merit function being not differentiable everywhere, the Eq.(S35) describes the directional derivative of $\phi_1\left(\mathbf{x}^k; \boldsymbol{\lambda}_{qp}^k, \boldsymbol{\mu}_{qp}^k\right)$ along the direction $\mathbf{d}^k$ generated by the SQP subproblem.

$$D\left(\phi_1\left(\mathbf{x}^k; \boldsymbol{\lambda}_{qp}^k, \boldsymbol{\mu}_{qp}^k\right); \mathbf{d}^k\right) = \nabla f\left(\mathbf{x}^k\right)^\top \mathbf{d}^k - \left(\boldsymbol{\rho}^k\right)^\top \left|\mathbf{h}\left(\mathbf{x}^k\right)\right| - \left(\mathbf{v}^k\right)^\top \mathbf{g}\left(\mathbf{x}^k\right)^- \tag{S35}$$

**Lemma 1.** *Let* $\mathbf{d}^k$, $\boldsymbol{\lambda}_{qp}^k$ *and* $\boldsymbol{\mu}_{qp}^k$ *be generated by problem* (**QP**). *Then the directional derivative of* $\phi_1\left(\mathbf{x}^k; \boldsymbol{\lambda}_{qp}^k, \boldsymbol{\mu}_{qp}^k\right)$ *in the direction* $\mathbf{d}^k$ *satisfies:*

$$D\left(\phi_1\left(\mathbf{x}^k; \boldsymbol{\lambda}_{qp}^k, \boldsymbol{\mu}_{qp}^k\right); \mathbf{d}^k\right) \leq 0 \tag{S36}$$

**Proof.** By applying Taylor's theorem to $f$, $\mathbf{h}$ and $\mathbf{g}$, we obtain:

$$
\begin{aligned}
\phi_1\left(\mathbf{x}^k + \alpha \mathbf{d}^k; \boldsymbol{\lambda}_{qp}^k, \boldsymbol{\mu}_{qp}^k\right) - \phi_1\left(\mathbf{x}^k; \boldsymbol{\lambda}_{qp}^k, \boldsymbol{\mu}_{qp}^k\right) \quad &= f\left(\mathbf{x}^k + \alpha \mathbf{d}^k\right) - f^k \\
&\quad + \left|\mathbf{h}\left(\mathbf{x}^k + \alpha \mathbf{d}^k\right)\right|^\top \boldsymbol{\rho}^k - \left|\mathbf{h}^k\right|^\top \boldsymbol{\rho}^k \\
&\quad + \left(\mathbf{v}^k\right)^\top \mathbf{g}\left(\mathbf{x}^k + \alpha \mathbf{d}^k\right)^- - \left(\mathbf{v}^k\right)^\top \mathbf{g}^{k-} \\
&\leq \alpha \nabla f^{k\top} \mathbf{d}^k + \gamma \alpha^2 \left\|\mathbf{d}^k\right\|_2^2 \\
&\quad + \left(\left|\mathbf{h}^k + \alpha \left(\nabla \mathbf{h}^k\right)^\top \mathbf{d}^k\right|^\top - \left|\mathbf{h}^k\right|^\top\right)\boldsymbol{\rho}^k \\
&\quad + \left(\max\left(0, -\mathbf{g}^k - \alpha \left(\nabla \mathbf{g}^k\right)^\top \mathbf{d}^k\right) - \left(\mathbf{g}^{k-}\right)^\top\right)\mathbf{v}^k
\end{aligned} \tag{S37}
$$

where $f^k := f(\mathbf{x}^k)$, $\mathbf{h}^k := \mathbf{h}(\mathbf{x}^k)$, $\mathbf{g}^{k-} := max\left(0, -\mathbf{g}(\mathbf{x}^k)\right)$, the positive constant $\gamma$ bounds the second-derivative terms in objective function and constraints, and $\alpha \in (0,1)$ is the step length.

According to Karush-Kuhn-Tucker (KKT) condition[16] of subproblem (QP), we have that $\mathbf{h}^k + \alpha(\nabla \mathbf{h}^k)^\top \mathbf{d}^k = 0$ and $\mathbf{g}^k + \alpha(\nabla \mathbf{g}^k)^\top \mathbf{d}^k \leq 0$, so for $\alpha \leq 1$ we have that

$$\phi_1\left(\mathbf{x}^k + \alpha\mathbf{d}^k; \boldsymbol{\lambda}_{qp}^k, \boldsymbol{\mu}_{qp}^k\right) - \phi_1\left(\mathbf{x}^k; \boldsymbol{\lambda}_{qp}^k, \boldsymbol{\mu}_{qp}^k\right) \begin{aligned}[t] &\leq \alpha\nabla f^{k\top}\mathbf{d}^k + \gamma\alpha^2\left\|\mathbf{d}^k\right\|_2^2 \\ &\quad + \left(\left(1-\alpha\right)\left|\mathbf{h}^k\right| - \left|\mathbf{h}^k\right|\right)^{\top}\boldsymbol{\rho}^k \\ &\quad + \left(\left(1-\alpha\right)\mathbf{g}^{k-} - \mathbf{g}^{k-}\right)^{\top}\mathbf{v}^k \\ &= \alpha\left[\nabla f^{k\top}\mathbf{d}^k - \left|\mathbf{h}^k\right|^{\top}\boldsymbol{\rho}^k - \left(\mathbf{g}^{k-}\right)^{\top}\mathbf{v}^k\right] \\ &\quad + \gamma\alpha^2\left\|\mathbf{d}^k\right\|_2^2 \end{aligned} \tag{S38}$$

Similarly, the following lower bound can be obtained:

$$\phi_1\left(\mathbf{x}^k + \alpha\mathbf{d}^k; \boldsymbol{\lambda}_{qp}^k, \boldsymbol{\mu}_{qp}^k\right) - \phi_1\left(\mathbf{x}^k; \boldsymbol{\lambda}_{qp}^k, \boldsymbol{\mu}_{qp}^k\right) \geq \alpha\left[\nabla f^{k\top}\mathbf{d}^k - \left|\mathbf{h}^k\right|^{\top}\boldsymbol{\rho}^k - \left(\mathbf{g}^{k-}\right)^{\top}\mathbf{v}^k\right] - \gamma\alpha^2\left\|\mathbf{d}^k\right\|_2^2 \tag{S39}$$

The directional derivative of $\phi_1$ in the direction $\mathbf{d}^k$ is given by Eq.(S32) by taking the limit. Due to $\mathbf{d}^k$ satisfying the KKT condition of the subproblem (**QP**), we can reformulate the directional derivative as:

$$D\left(\phi_1\left(\mathbf{x}^k; \boldsymbol{\lambda}_{qp}^k, \boldsymbol{\mu}_{qp}^k\right); \mathbf{d}^k\right) \begin{aligned}[t] &= -\left(\mathbf{d}^k\right)^{\top}\cdot\mathbf{B}^k\cdot\mathbf{d}^k + \left(\mathbf{d}^k\right)^{\top}\cdot\left(\nabla\mathbf{h}^k\right)^{\top}\cdot\boldsymbol{\lambda}_{qp}^{k+1} + \mathbf{d}^{k\top}\nabla\mathbf{g}^{k\top}\boldsymbol{\mu}_{qp}^{k+1} \\ &\quad - \left(\boldsymbol{\rho}^k\right)^{\top}\left|\mathbf{h}\left(\mathbf{x}^k\right)\right| - \left(\mathbf{v}^k\right)^{\top}\mathbf{g}\left(\mathbf{x}^k\right)^{-} \\ &\leq -\mathbf{d}^{k\top}\mathbf{B}^k\mathbf{d}^k + \left|\mathbf{h}^k\right|^{\top}\left(\boldsymbol{\lambda}_{qp}^{k+1} - \boldsymbol{\rho}^k\right) + \left(\mathbf{g}^{k-}\right)^{\top}\left(\boldsymbol{\mu}_{qp}^{k+1} - \mathbf{v}^k\right) \\ &\leq -\mathbf{d}^{k\top}\mathbf{B}^k\mathbf{d}^k + \left|\mathbf{h}^k\right|^{\top}\left(\boldsymbol{\lambda}_{qp}^{k+1} - \boldsymbol{\rho}^k\right) + \left(\mathbf{g}^{k-}\right)^{\top}\left(\boldsymbol{\mu}_{qp}^{k+1} - \mathbf{v}^k\right) \end{aligned} \tag{S40}$$

Since $\mathbf{B}^k$ is positive definite, and according to Eq.(S33) and Eq.(S34), we have $\boldsymbol{\rho}^k \geq \boldsymbol{\lambda}_{qp}^{k+1}$ and $\mathbf{v}^k \geq \boldsymbol{\mu}_{qp}^{k+1}$, then the Lemma 1 has been proven. □

Lemma 1 indicates that if $\mathbf{d}^k$ is the KKT point of problem (**QP**) and $\mathbf{B}^k$ is positive definite, then the merit function decreases. This implies that in each iteration, either the infeasibility of the original problem is alleviated, or the value of the objective function is reduced, or both situations occur simultaneously.

The Armijo condition guarantees a sufficient decrease of the merit function:

$$\phi_1\left(\mathbf{x}^k + \alpha\mathbf{d}^k; \boldsymbol{\lambda}_{qp}^k, \boldsymbol{\mu}_{qp}^k\right) - \phi_1\left(\mathbf{x}^k; \boldsymbol{\lambda}_{qp}^k, \boldsymbol{\mu}_{qp}^k\right) < \alpha\eta D\left(\phi_1\left(\mathbf{x}^k; \boldsymbol{\lambda}_{qp}^k, \boldsymbol{\mu}_{qp}^k\right); \mathbf{d}^k\right) \tag{S41}$$

Where $\eta \in (0, 0.5)$ is a hyper-parameter. However, it is possible that no step length $\alpha$ can satisfy the Armijo condition, regardless of how small it is. Under such circumstances, we accept the last $\alpha$ when the maximum number of line searches has been reached.

For practical engineering problems, due to numerical stability and other factors, the KKT conditions as criteria may be too strict[17], making the algorithm difficult to achieve convergence. Therefore, based on previous research[15,18,19], we check whether the following criteria are met after every line search:

$$\left\| \mathbf{h}\left( \mathbf{x}^k + \alpha \mathbf{d}^k \right) \right\|_1 + \left\| \mathbf{g}\left( \mathbf{x}^k + \alpha \mathbf{d}^k \right)^- \right\|_1 < tol \tag{S42}$$

$$\left| f\left( \mathbf{x}^k + \alpha \mathbf{d}^k \right) - f\left( \mathbf{x}^k \right) \right| < tol \tag{S43}$$

$$\left\| \mathbf{d}^k \right\|_2 < tol \tag{S44}$$

where, $tol$ represents the convergence tolerance. The solution is considered optimal if Eqs.(S42), (S43) or Eqs.(S42), (S44) are satisfied. The satisfaction of Eq.(S42) indicates that the optimal solution is feasible. The satisfaction of either Eq.(S43) or Eq.(S44) suggests that there is little potential for further decrease in the objective function.

The SQP-driven feasible path algorithm for surrogate-assisted BBO is as shown in ALGORITHM S2.

ALGORITHM S2 SQP-Driven feasible path algorithm for surrogate-assisted BBO.

---

**Data:** Surrogate model $s(\cdot)$, $\boldsymbol{\lambda}_{qp}^0 \leftarrow 0, \boldsymbol{\mu}_{qp}^0 \leftarrow 0, \boldsymbol{\rho}^0 \leftarrow 0, \boldsymbol{\nu}^0 \leftarrow 0,$
$\eta \leftarrow 0.1, tol \leftarrow 10^{-6}$, line search multiple $\beta \leftarrow 0.618$, line search
maximum iteration $n_{max} \leftarrow 10$

**Input:** feasible initial point $\boldsymbol{x}^0$

**Output:** solution $\boldsymbol{x}^*$

1 **begin**
2    **while** $acc_{inf} > tol \vee (acc_{step} > tol \wedge acc_{opt} > tol)$ **do**
3      $k \leftarrow k + 1$
4      $\boldsymbol{y}^k \leftarrow \boldsymbol{s}(\boldsymbol{x}^k)$
5      Calculate $\nabla_x \boldsymbol{y}^k$ and $\nabla_x^2 \boldsymbol{y}^k$ based on Section 3
6      Calculate $\nabla_x f$, $\nabla_x \boldsymbol{g}^k$, $\nabla_x \boldsymbol{h}^k$ and $\boldsymbol{H}^k$
7      **if** $\boldsymbol{H}^k$ *is positive definite* **then**
8        $\boldsymbol{B}^k \leftarrow \boldsymbol{H}^k$
9      **else**
10        $\boldsymbol{B}^k \leftarrow \boldsymbol{H}^k + \boldsymbol{E}^k$
11      **end if**
12      Solve (QP) subproblem, then obtain $\boldsymbol{d}^k$, $\boldsymbol{\lambda}_{qp}^k$ and $\boldsymbol{\mu}_{qp}^k$
13      $\boldsymbol{\rho}^k \leftarrow max(|\boldsymbol{\lambda}_{qp}^k|, \frac{|\boldsymbol{\lambda}_{qp}^k| + \boldsymbol{\rho}^{k-1}}{2})$
14      $\boldsymbol{\nu}^k \leftarrow max(|\boldsymbol{\mu}_{qp}^k|, \frac{|\boldsymbol{\mu}_{qp}^k| + \boldsymbol{\nu}^{k-1}}{2})$
15      Calculate $\phi_1(\boldsymbol{x}^k, \boldsymbol{\lambda}_{qp}^k, \boldsymbol{\mu}_{qp}^k)$ and $D(\phi_1(\boldsymbol{x}^k, \boldsymbol{\lambda}_{qp}^k, \boldsymbol{\mu}_{qp}^k); \boldsymbol{d}^k)$
16      $n \leftarrow 0$, $\alpha \leftarrow 1$
17      **while** $n \leq n_{max}$ **do**
18        $n \leftarrow n + 1$
19        $\boldsymbol{x}^{k+1} \leftarrow \boldsymbol{x}^k + \alpha \boldsymbol{d}^k$
20        **if** $\phi_1(\boldsymbol{x}^k + \alpha\boldsymbol{d}^k, \boldsymbol{\lambda}_{qp}^k, \boldsymbol{\mu}_{qp}^k) - \phi_1(\boldsymbol{x}^k, \boldsymbol{\lambda}_{qp}^k, \boldsymbol{\mu}_{qp}^k) <$
$\alpha\eta D(\phi_1(\boldsymbol{x}^k, \boldsymbol{\lambda}_{qp}^k, \boldsymbol{\mu}_{qp}^k); \boldsymbol{d}^k)$ **then**
21          **break**
22        **else**
23          $\alpha \leftarrow \beta\alpha$
24        **end if**
25      **end while**
26    **end while**
27 **end**

## S4. Examples 1-6

The sampling range of the six test functions in Examples 1-6 is selected as [–2, 2], with a training set sample size of 2000, a validation set of 50, and a test set of 200.

TABLE S2 Test functions for Examples 1-6

| Example | Dimension | Expression | Minimal value | Minimum point |
|---|---|---|---|---|
| 1 | 10 | $$\sum_{i=1}^{10} x_i^2$$ | 0 | $(0,0,\ldots,0)$ |
| 2 | 10 | $$\sum_{i=1}^{10} x_i^2 + \sum_{i=1}^{9}(x_i - x_{i+1})^2$$ | 0 | $(0,0,\ldots,0)$ |
| 3 | 2 | $$\left(4 - 2.1x_1^2 + \frac{x_1^4}{3}\right)x_1^2 + x_1 x_2 + 4(-1 + x_2^2)x_2^2$$ | –1.03 | $(0.09,–0.71)$ $(–0.09,0.71)$ |
| 4 | 2 | $$0.5 + \frac{\sin^2(x_1^2 - x_2^2) - 0.5}{[1 + 0.001(x_1^2 + x_2^2)]^2}$$ | 0 | $(0,0)$ |
| 5 | 5 | $$1 + \frac{1}{4000}\sum_{i=1}^{5} x_i^2 - \prod_{i=1}^{5} \cos\left(\frac{x_i}{\sqrt{i}}\right)$$ | 0 | $(0,0,0,0,0)$ |
| 6 | 5 | $$-20\,exp\left(-0.2\sqrt{\frac{1}{5}\sum_{i=1}^{5} x_i^2}\right) - exp\left(\frac{1}{5}\sum_{i=1}^{5} \cos(2\pi x_i)\right) + 20 + e$$ | 0 | $(0,0,0,0,0)$ |

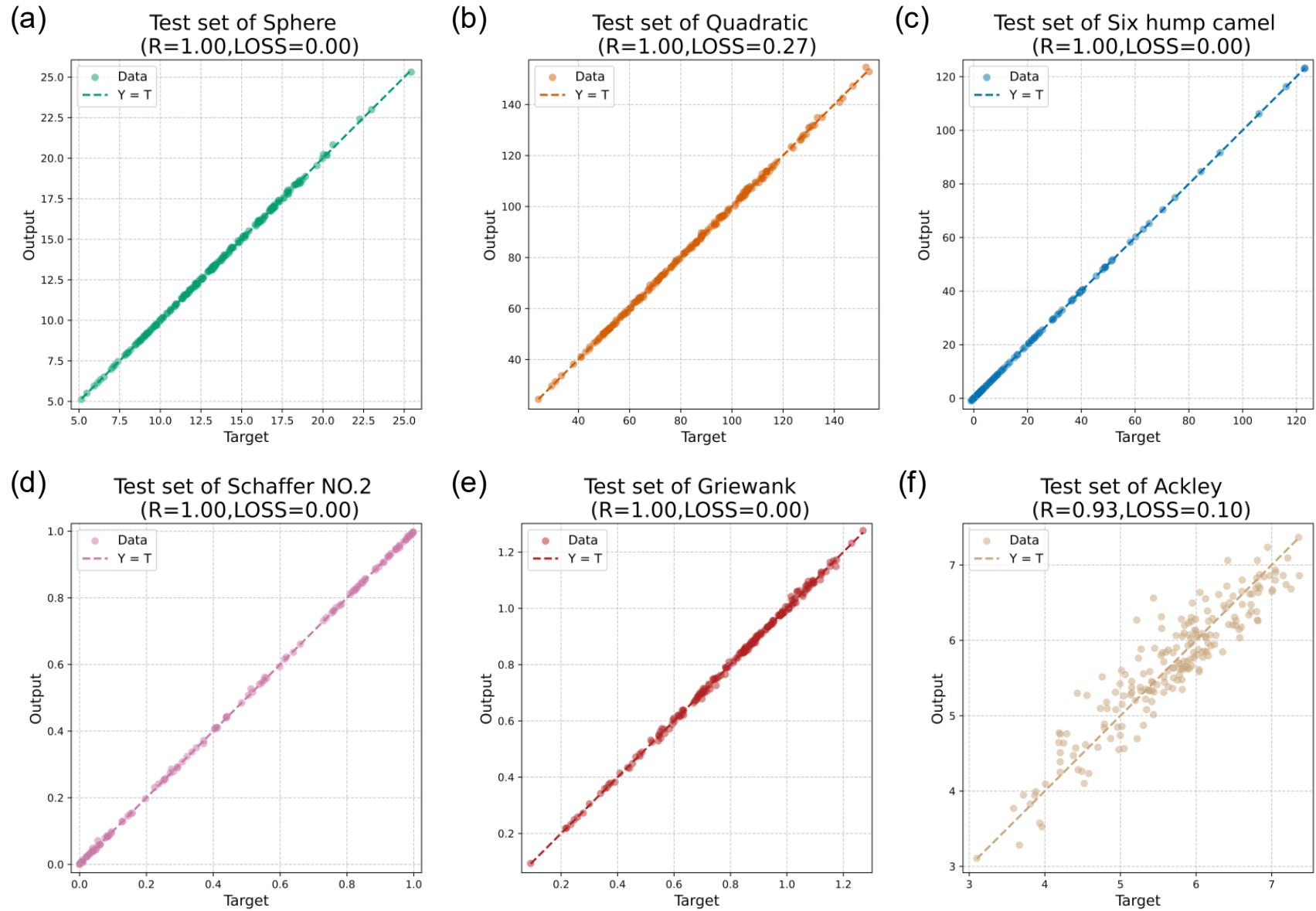

FIGURE S17 Parity plots of (a) Sphere function, (b) Quadratic function, (c) Six hump camel function, (d) Schaffer function NO.2, (e) Griewank function and (f) Ackley function on the test set

## S5. Examples 7-10

In these examples, we employ 4 classic test functions as black-box, and compare the performance of MLFP and TRF on these functions.

TABLE S3 Test functions for Examples 7-10

| Example | Dimension | Expression | Minimal value | Minimum point |
|---------|-----------|------------|---------------|---------------|
| 7 | 3 | $-\sum_{i=1}^{4}\alpha_i exp\left(-\sum_{j=1}^{3}A_{ij}\left(x_j - P_{ij}\right)^2\right)$ | −3.8628 | (0.1146, 0.5556, 0.8525) |
| 8 | 4 | $(x_1+10x_2)^2+5(x_3-x_4)^2$ $+(x_2-2x_3)^4+10(x_1-x_4)^4$ | 0 | (0,0,…,0) |
| 9 | 4 | $\sum_{i=1}^{4}[100(x_{i+1}-x_i^2)^2+(x_i-1)^2]$ | 0 | (1,1,…,1) |
| 10 | 6 | $\sum_{i=1}^{6}(x_i-1)^2-\sum_{i=2}^{6}x_ix_{i-1}$ | −50 | $x_i$ $= i(7-i)$ $for\ i$ $= 1,…,6$ |

Note: $\alpha = [1\quad 1.2\quad 3\quad 3.2]^T$, $A = \begin{bmatrix} 3 & 10 & 30 \\ 0.1 & 10 & 35 \\ 3 & 10 & 30 \\ 0.1 & 10 & 35 \end{bmatrix}$,

$P = 10^{-4}\begin{bmatrix} 3689 & 1170 & 2673 \\ 4699 & 4387 & 7470 \\ 1091 & 8732 & 5547 \\ 381 & 5743 & 8828 \end{bmatrix}$

## S6. Example 11: Williams-Otto process

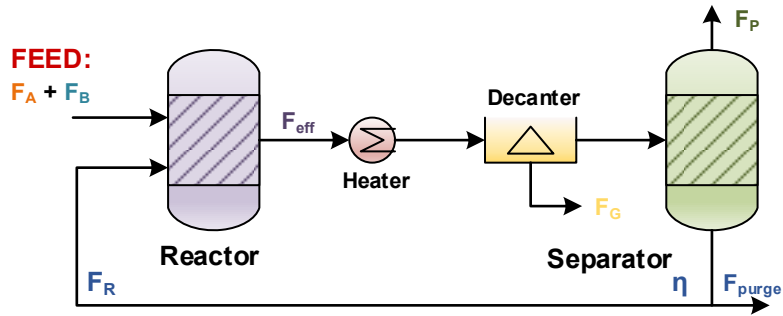

FIGURE 18 Diagram of Williams-Otto process.

The Williams-Otto (WO) process within the reactor involves the following 3 reactions[20,21]:

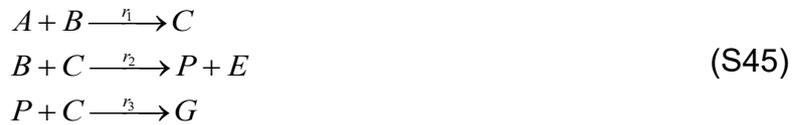

$$A + B \xrightarrow{\ r_1\ } C$$
$$B + C \xrightarrow{\ r_2\ } P + E$$
$$P + C \xrightarrow{\ r_3\ } G$$

(S45)

All reactants involved are assumed to be fictitious components, and the mixture is considered to have a constant density of $\rho = 50$. The separation efficiency is also assumed to be fixed. Reactants A and B are introduced into the reactor to produce P with by-product C, E, and G. This stream is then sent to a decanter, where component G is removed. The remaining mixture is fed into a distillation column, with the desired product P obtained as the top distillate. A portion of the bottom stream is recycled back to the reactor to improve process efficiency.

The full equation-oriented optimization formulation is presented as follows:

**Objective:**

$$\max ROI = \frac{2207 F_P + 50 F_{purge} - 168 F_A - 252 F_B - 2.22 F_{eff}^{sum} - 84 F_G - 60 V \rho}{600 V \rho} \times 100\%$$

(S46)

**Kinetics:**

$$r_1 = 5.9755 \times 10^9 \exp\left(-120/T\right) x_A x_B V \rho$$
$$r_2 = 2.5962 \times 10^{12} \exp\left(-150/T\right) x_B x_C V \rho$$
$$r_3 = 9.6283 \times 10^{15} \exp\left(-200/T\right) x_P x_C V \rho$$

(S47)

**Reactor balance:**

$$F_{eff}^A = F_A + F_R^A - r_1$$
$$F_{eff}^B = F_B + F_R^B - (r_1 + r_2)$$
$$F_{eff}^C = F_R^C + 2r_1 - 2r_2 - r_3$$
$$F_{eff}^E = F_R^E + 2r_2$$
$$F_{eff}^P = 0.1F_R^E + r_2 - 0.5r_3 \tag{S48}$$
$$F_{eff}^G = 1.5r_3$$
$$F_{eff}^{sum} = \sum_i F_{eff}^j \qquad j \in \{A,B,C,E,G,P\}$$
$$F_{eff}^j = F_{eff}^{sum} x_j \qquad j \in \{A,B,C,P\}$$

**Waste stream:**

$$F^G = F_{eff}^G \tag{S49}$$

**Product stream:**

$$F_P = F_{eff}^P - 0.1F_{eff}^E$$
$$F_{purge} = \eta\left(F_{eff}^A + F_{eff}^B + F_{eff}^C + 1.1F_{eff}^E\right) \tag{S50}$$

**Recycle stream:**

$$F_R^j = (1-\eta)F_{eff}^j \qquad j \in \{A,B,C,E\}$$
$$F_R^P = 0.1(1-\eta)F_{eff}^E \tag{S51}$$

**Bound constraints:**

$$V \in [0.03,0.1] \quad T \in [5.8,6.8] \quad F_P \in [0,4.763]$$
$$F \geq 0 \qquad\qquad F_A, F_B \geq 1 \qquad \eta \in [0,1] \tag{S52}$$

Eq.(S46) -(S52) is the EO-based WO optimization problem. Establishing a surrogate model for Eq.(S46)-(S51) with the input of $x = [V, T, \eta, F_A, F_B]$ and an the output of $y = [-ROI, F_P]$, it can be expressed as the following optimization problem:

$$\begin{aligned}&\min_{\mathbf{x}\in\mathbb{R}^5} \quad -ROI \\ &\text{s.t.} \quad \mathbf{y} = s(\mathbf{x}) \\ &\qquad\quad \text{Eq.(S41)}\end{aligned} \tag{S53}$$

## S7. Example12: Extractive distillation for toluene and n-heptane separation

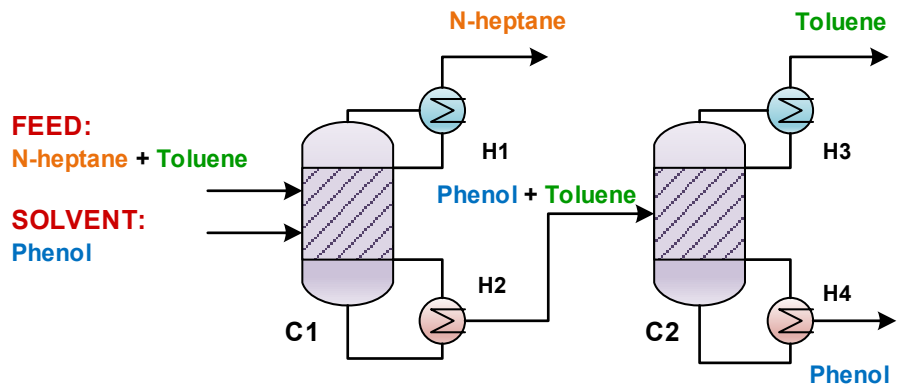

FIGURE S19 Process flow diagram for separation of toluene form n-heptane using solvent phenol in Ma et al.'s work[22].

In this example, the sampling numbers for the training set, validation set, and test set were set to 2000, 200, and 1000, with 1014, 95, and 511 successful points, respectively. The multilayer perceptron (MLP) is employed as surrogate model in this case, who possesses four hidden layers with hidden layers 1-2 having 100 neurons and hidden layers 3-4 having 120 neurons. FIGURE S20 shows the performance of constructed model on test set. As is shown in FIGURE S8, the surrogate model can predict the output perfectly, since each graph exhibits a R of 1.00 and relative errors are no greater than 0.56%. Such results are highly encouraging for model validation and instill strong confidence in the models' reliability.

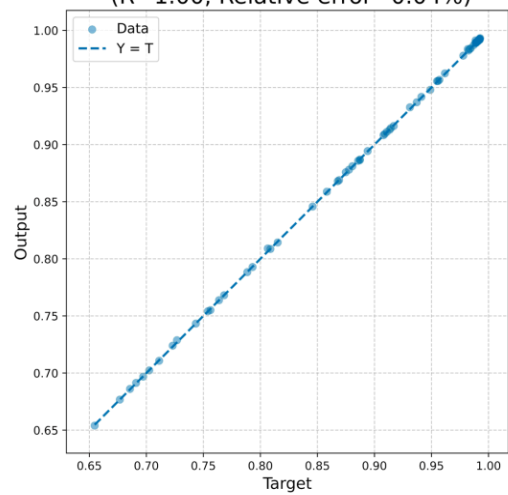 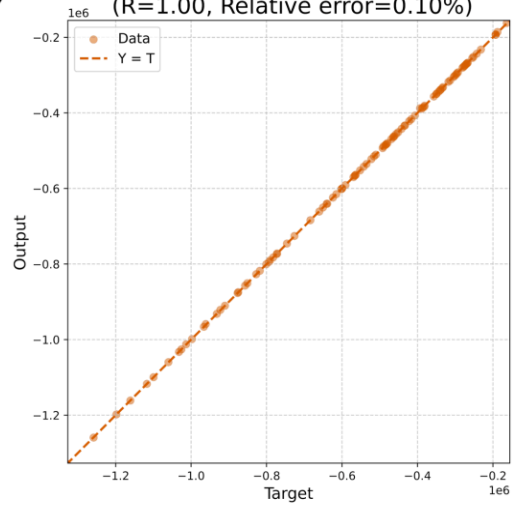 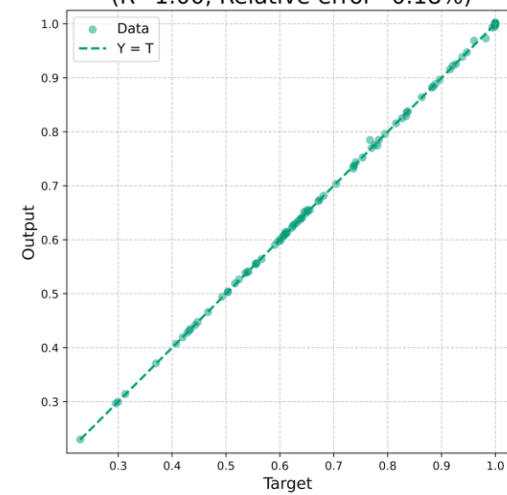

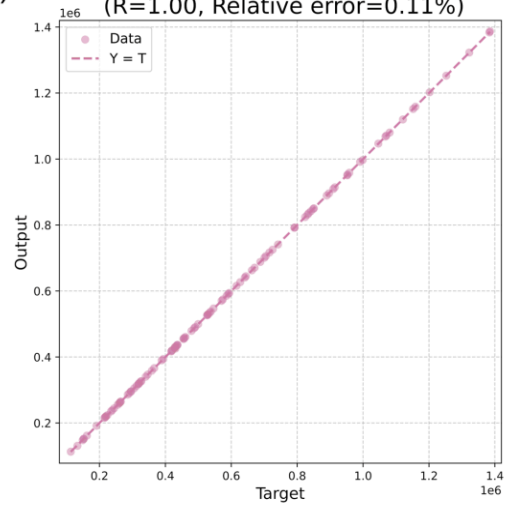 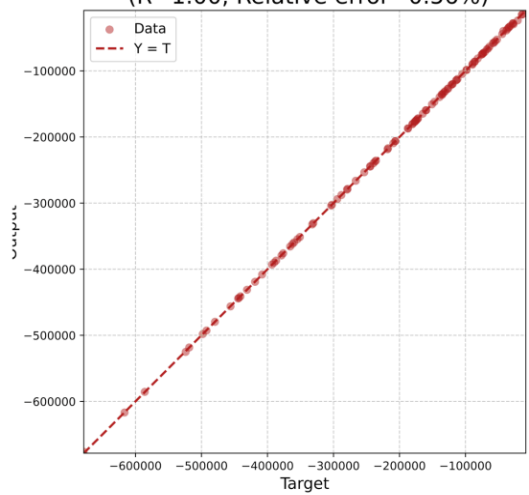 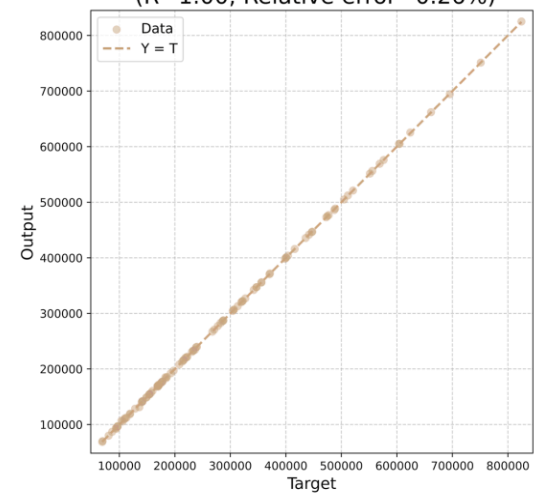

FIGURE S20 Parity plots of (a) mole fraction of n-heptane in C1 distillate, (b) mole fraction of toluene in C2 distillate, (c) condenser duty of C1, (d) reboiler duty of C1, (e) condenser duty of C2, and (f) reboiler duty of C2 on the test set

## S8. Example 13: CO₂ capture

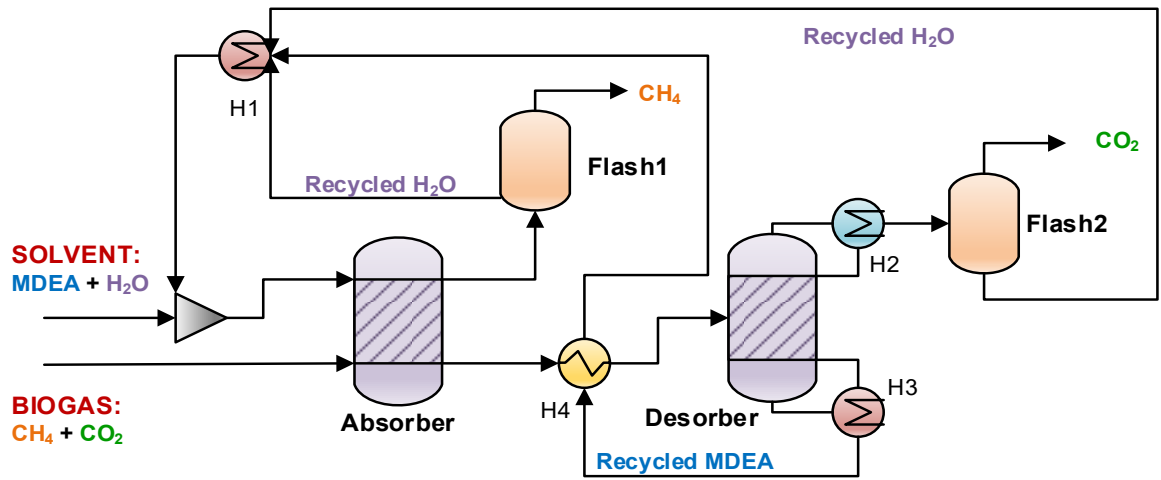

FIGURE S21 Process flow diagram of CO₂ capture from biogas using aqueous

MDEA

The surrogate model for the entire process is also a MLP, possessing six hidden layers with 1-3 layers having 100 neural and 4-6 layer having 200 neural. The sampling numbers for the training set, validation set, and test set were set to 2000, 100, and 200.

The performance on test set as is shown in FIGURE S22, where each graph exhibits a R of 1.00 and relative errors are no greater than 0.10%.

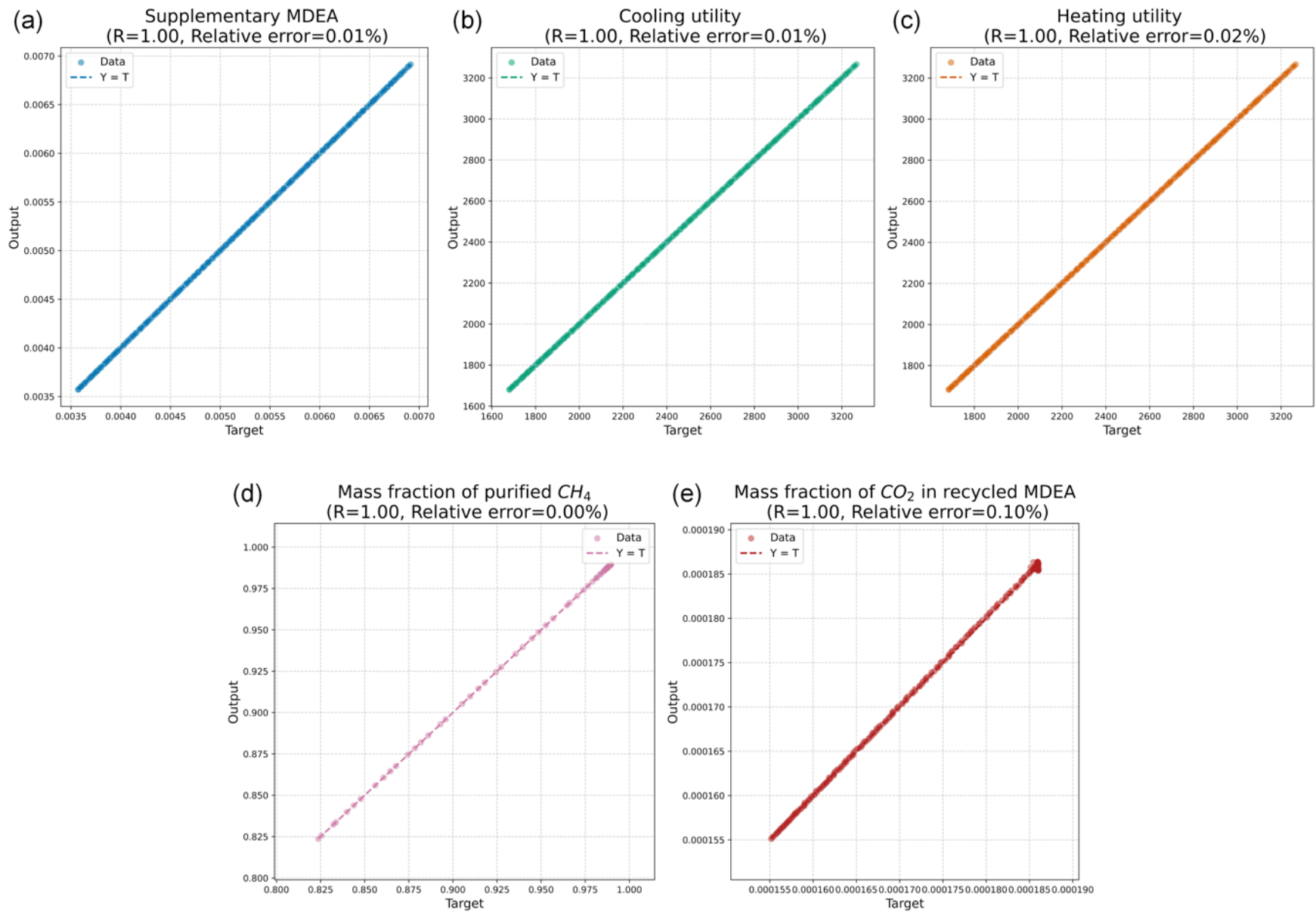

FIGURE S22 Parity plots of (a) supplementary MDEA, (b) cooling utility, (c) heating utility, (d) mass fraction of CH₄ in the purified

biogas (e) mass fraction of $CO_2$ in the mean MEDA, and (f) reboiler duty of column C2